\documentclass{svproc}
\usepackage{url}

\usepackage{graphicx}
\usepackage{multicol}
\usepackage{footmisc}

\usepackage{amsfonts}
\usepackage{amssymb}
\usepackage{amsmath}
\usepackage{tikz}
\usepackage{listings} 
\usepackage{centernot} 
\usepackage{hyperref} 
\usepackage{varwidth} 
\usepackage[misc]{ifsym} 
\usepackage{multirow} 

\usepackage{array,etoolbox}
\preto\tabular{\setcounter{magicrownumbers}{0}}
\newcounter{magicrownumbers}
\def\rownumber{}
\preto\tabular{\setcounter{magicrownumbers2}{0}}
\newcounter{magicrownumbers2}
\def\foo{}

\newcommand{\itparen}[1]{\textup{(}#1\textup{)}}

\newcommand{\Distr}[1]{\mathrm{Distr} \: #1}
\newcommand{\mDistr}[1]{\mathrm{mDistr} \: #1}
\newcommand{\jDistr}[1]{\mathrm{jDistr} \: #1}

\spnewtheorem{algorithm}{Algorithm}{}{}

\begin{document}
\mainmatter              
\title{$ST$-distributive and $ST$-modular Lattices}
\titlerunning{$ST$-modular and $ST$-distributive Lattices}  
%
\author{M. R. Emamy-K. \textsuperscript{[0000-0002-5407-0392]} \and \\ Gustavo A. Mel\'{e}ndez R\'{i}os \Letter \hspace{.3mm} \textsuperscript{[0000-0001-5980-3666]}}
\authorrunning{M. R. Emamy-K. \and Gustavo A. Mel\'{e}ndez R\'{i}os} 
%
\tocauthor{M. R. Emamy-K., Gustavo A. Mel\'{e}ndez R\'{i}os}
\institute{Universidad de Puerto Rico, Recinto de R\'{i}o Piedras, San Juan PR, USA, \\
\email{mreza.emamy@upr.edu, gustavo.melendez@upr.edu},\\ 
}

\maketitle              

\begin{abstract}
For two subsets $S$ and $T$ of a given lattice $L,$ we define 
a relative distributive (modular) property over $L,$ that underlies a  
large family including the usual class of distributive (modular) lattices. Our  proposed class will be called $ST$-distributive ($ST$-modular) lattices. In this paper, we find examples of maximal $S$ and $T$ to form $ST$-distributive lattices for several non-distributive finite lattices of small order. In particular, we characterize the maximal pairs of subsets $(S,T)$, subject to certain constraints, that induce $ST$-distributivity in the lattice family $\mathbf{M}_{n,n}$. Afterwards, we present an application of $ST$-modular to convex sets and polytopes. This application is the original motivation for our new definitions. The aforementioned definitions are closely related to distributive elements and Stanley's SS-lattices.
\end{abstract}

\section{Introduction}

We propose two new classes of lattices: $ST$-distributive lattices and $ST$-modular lattices. The idea behind them is to define relative distributive and modular properties that are fulfilled by certain pairs of subsets of a lattice that is not necessarily distributive nor modular. Naturally, we desire that these new classes include distributive and modular lattices respectively (i.e. are a generalization of them). In a sense, this is analogous to relative topologies in topological spaces.\\

\indent This paper introduces one of the most basic questions regarding these new classes: Given a lattice, which pairs of subsets satisfy these relative distributive/modular properties? Of course, this question is non-trivial only if the given lattice is non-distributive/non-modular. The main goal coming out of this question is to identify families of lattices with particular characterizations of the subsets that induce these relative properties. \\

\indent Here, we define these new classes precisely and present some initial results on $ST$-distributivity for a family of lattices called $\mathbf{M}_{n,n}$. We also show connections between $ST$-distributive lattices and two established concepts in lattice theory: distributive elements \cite{G2011} and Stanley's SS-lattices \cite{S}. In addition, we present a potential applications of $ST$-modular lattices to convex sets. \\

\indent This paper is organized as follows: Section \ref{sec_st-distr} discusses $ST$-distributive lattices including our results on $\mathbf{M}_{n,n}$. Afterwards, Sect. \ref{sec_st-mod} introduces $ST$-modular lattices and applies them to convex sets. Finally, Sect. \ref{sec_concl} concludes with a summary and ideas for continuing this work. For more details on how Sects. \ref{sec_st-distr} and \ref{sec_st-mod} are organized internally, see their respective introductory paragraphs.

\section{$ST$-distributive Lattices} \label{sec_st-distr}

In this section, we discuss $ST$-distributive lattices. Subsection \ref{basic-def-st-distr} establishes the basic definitions regarding this new concept and provides some illustrative examples. This is developed further with some basic properties of these lattices in Subsect. \ref{basic_props}. Subsection \ref{search} then specifies our search problem and presents the methodology we use to tackle it. Afterwards, the family of lattices to which we apply this methodology, $\mathbf{M}_{n,n}$, is described in Subsect. \ref{fam-M_n,n}. This is followed by the results of this application in Subsect. \ref{results}. Finally, Subsections \ref{distr-elem-connect} and \ref{ss-lat-connect} conclude this discussion by connecting our $ST$-distributive lattices to distributive elements and Stanley's SS-lattices respectively.

\subsection{Basic Definitions} \label{basic-def-st-distr}

We introduce $ST$-distributive lattices. A lattice $L$ is called $ST$-distributive if it satisfies a relative distributive property with respect to two subsets $S,T \subseteq L$. The idea is that all elements of $S$ can be distributed into any two elements of $T$ in both ways: meet into join and join into meet. We now provide a formal definition.

\begin{definition} ($ST$-distributive Lattice). \label{st-distr}
Given a lattice $L$ with $S, T \subseteq L$, we define:

\begin{itemize}
\item \underline{$ST$-meet Distributive Lattice}:
$L$ is said to be $ST$-meet distributive if for all  $s \in S$ and $t_1, t_2 \in T$,
\begin{equation} \label{st-meet_distr}
s \land \left( t_1 \lor t_2 \right) = \left( s \land t_1 \right) \lor \left( s \land t_2 \right) \; .
\end{equation}

\item \underline{$ST$-join Distributive Lattice}:
$L$ is said to be $ST$-join distributive if for all  $s \in S$ and $t_1, t_2 \in T$,
\begin{equation} \label{st-join_distr}
s \lor \left( t_1 \land t_2 \right) = \left( s \lor t_1 \right) \land \left( s \lor t_2 \right) \; .
\end{equation}

\item \underline{$ST$-distributive Lattice}:
$L$ is said to be $ST$-distributive if it is both $ST$-meet distributive  and $ST$-join distributive.
\end{itemize}
\end{definition}


\begin{figure}
\begin{minipage}{0.45\textwidth}
\centering
\begin{tikzpicture}
[roundnode/.style={circle, draw=black, fill=white, very thin, text width =4.00mm, inner sep = 0pt, text centered}]
 \node [roundnode] (1) at (0,2) {$1$};
  \node [roundnode] (u) at (-1,1.4) {$u$};
  \node [roundnode] (v) at (-1,.6) {$v$};
  \node [roundnode] (w) at (1,1) {$w$};
  \node [roundnode] (0) at (0,0) {$0$};
 \draw (1) -- (u) -- (v) -- (0) -- (w) -- (1);
\end{tikzpicture}
\caption{The pentagon $\mathbf{N}_5$} \label{n5}
\end{minipage}
\begin{minipage}{0.45\textwidth}
\centering
\begin{tikzpicture}
[roundnode/.style={circle, draw=black, fill=white, very thin, text width =4.00mm, inner sep = 0pt, text centered}]
 \node [roundnode] (1) at (0,2) {$1$};
  \node [roundnode] (a) at (-1,1) {$a$};
  \node [roundnode] (b) at (0,1) {$b$};
  \node [roundnode] (c) at (1,1) {$c$};
  \node [roundnode] (0) at (0,0) {$0$};
 \draw (1) -- (a) -- (0) -- (b) -- (1) -- (c) -- (0);
\end{tikzpicture}
\caption{The diamond $\mathbf{M}_3$} \label{m3}
\end{minipage}
\end{figure}

\begin{figure}
\centering
\begin{tikzpicture}
[roundnode/.style={circle, draw=black, fill=white, very thin, text width =4.00mm, inner sep = 0pt, text centered}]
 \node [roundnode] (1) at (0,3.5) {$1$};
  \node [roundnode] (a) at (0,2) {$a$};
  \node [roundnode] (b) at (-2,2) {$b$};
  \node [roundnode] (c) at (2,2) {$c$};
  \node [roundnode] (d) at (-1,1) {$d$};
  \node [roundnode] (e) at (1,1) {$e$};
  \node [roundnode] (0) at (0,0) {$0$};
 \draw (1) -- (b) -- (d) -- (0) -- (e) -- (c) -- (1) -- (a) -- (d);
 \draw (a) -- (e);
\end{tikzpicture}
\caption{Lattice $L$ of Example \ref{meet-join}} \label{counter-ex-lattice}
\end{figure}
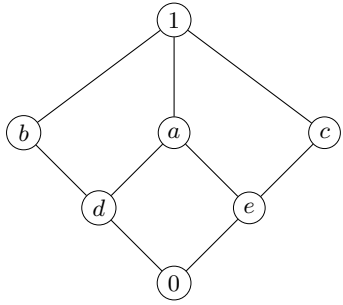

A couple of remarks are in order. First, observe that this definition generalizes distributive lattices which happen to be lattices that are $ST$-distributive for $S = T = L$. Second, every lattice $L$ is vacuously $ST$-distributive given $S = \emptyset$ or $T = \emptyset$. Third, we will occasionally use the term relative distributive property to refer to the $ST$-distributivity of a lattice. We now give an example of an $ST$-distributive lattice.

\begin{example} (Example of $ST$-distributive Lattice). \label{st-distr-ex}
Consider the lattice $\mathbf{N}_5$ shown in Fig. \ref{n5}. It is known to be non-distributive \cite{Bi,DP,G2011}. However, if $S = \left \{ u,v \right \}$ and $T = \left \{w,0 \right \}$, then it is $ST$-distributive. To verify this, we must check that it is both $ST$-meet distributive and $ST$-join distributive. The necessary computations are provided in Table \ref{n5_example_st-distr} for the case $s = u$ in Equations (\ref{st-meet_distr}) and (\ref{st-join_distr}). The computations for $s = v$ are identical.
\end{example}

\begin{table}
\caption{Computations to show $\mathbf{N}_5$ is $ST$-distributive in Example \ref{st-distr-ex}} \label{n5_example_st-distr}
\begin{center}
\begin{tabular}{ c m{7mm} c }
\hline
$ST$-meet & & $ST$-join \\
\hline 
$u \land \left( w \lor 0\right) = 0$ & &  $u \lor \left(w \land 0\right) = u$ \\
$\left( u \land w \right) \lor \left( u \land 0 \right) = 0$ & & $\left(u \lor w\right) \land \left(u \lor 0\right) = u$  \\
\hline
$u \land \left( 0 \lor w\right) = 0$ & & $u \lor \left(0 \land w\right) = u$ \\
$\left( u \land 0\right) \lor \left( u \land w\right) = 0$ & & $\left(u \lor 0\right) \land \left(u \lor w\right) = u$ \\
\hline
$u \land \left( 0 \lor 0\right) = 0$ & & $u \lor \left(0 \land 0\right) = u$ \\
$\left( u \land 0\right) \lor \left( u \land 0\right) = 0$ & & $\left(u \lor 0\right) \land \left(u \lor 0\right) = u$ \\
\hline
$u \land \left( w \lor w\right) = 0$ & & $u \lor \left(w \land w\right) = 1$ \\
$\left( u \land w\right) \lor \left( u \land w\right) = 0$ & & $\left(u \lor w\right) \land \left(u \lor w\right) = 1$\\
\hline
\end{tabular}
\end{center}
\end{table} 

\indent In practice, it will not be necessary to do all the computations done in Example \ref{st-distr-ex}. We will show that the order structure of the lattice will guarantee that certain triples $(s,t_1,t_2)$ will satisfy Equations (\ref{st-meet_distr}) and (\ref{st-join_distr}). See Property \ref{eff_crit} in Subsect. \ref{basic_props} for details.\\

\indent Next, we briefly discuss our new definition. We begin with a counter-example to show that $ST$-meet distributive and $ST$-join distributive are not equivalent in general. 

\begin{example} ($ST$-meet Distributive $\centernot\iff$ $ST$-join Distributive). \label{meet-join}
Consider the lattice $L =\left \{0, a, b, c, d, e, 1 \right \}$ shown in Fig. \ref{counter-ex-lattice}. Let $S =\left  \{a \right \}$ and $T=\left \{b, c \right \}$. Then $L$ is $ST$-meet distributive but not $ST$-join distributive.

\begin{proof}
\noindent \underline{$ST$-meet Distributive}:
There are four cases: 
\begin{equation}
a \land \left( b \lor c \right) \qquad a \land \left( c \lor b \right) \qquad a \land \left( b \lor b \right) \qquad a \land \left( c \lor c \right).
\end{equation}
Note that the second follows from the first by commutativity of lattice operations while the last two are immediate by idempotency. Hence, we need to only verify the first:
\begin{equation}
a \land \left(b \lor c\right) = a \land 1 = a = d \lor e = \left(a \land b\right) \lor \left(a \land c\right) \; .
\end{equation}

\noindent \underline{Not $ST$-join Distributive}: 
Note that
\begin{equation}
a \lor \left(b \land c\right) = a \lor 0 = a \neq 1 = 1 \land 1 = \left(a \lor b\right) \land \left(a \lor c\right) \; .
\end{equation}
\qed
\end{proof}
\end{example}

Another point to emphasize is that the subsets $S$ and $T$ in our definition are not interchangeable in general, or equivalently, that the pair $\left( S,T \right)$ is ordered. The following example illustrates this. 

\begin{example} \label{st-neq-ts} ($ST$-distributive $\centernot\iff$ $TS$-distributive).
Consider the lattice $\mathbf{M}_3$ shown in Fig. \ref{m3}. Let $S =\left \{a,b\right \}$ and $T = \left \{c \right \}$. Then $\mathbf{M}_3$ is $ST$-distributive but not $TS$-distributive.
\begin{proof}
\underline{ST-distributive}:
This follows from idempotency of lattice operations because $T$ only has one element. For instance, 
\begin{equation}
a \land \left(c \lor c\right) = a \land c = \left(a \land c\right) \lor \left(a \land c\right)  \; .
\end{equation}
\underline{Not TS-distributive}:
Observe that 
\begin{equation}
c \land \left(a \lor b\right) = c \land 1 = c \neq 0 = 0 \lor 0 = \left(c \land a\right) \lor \left(c \land b\right) \; .
\end{equation}
\qed
\end{proof}
\end{example}

We now introduce a special pair $\left( S,T \right)$ of proper subsets of a lattice: those that make it $ST$-distributive but that cannot be expanded further without either, breaking our relative distributive property or becoming the whole lattice. We call them maximal $ST$-pairs.

\begin{definition} \label{max_stp}(Maximal $ST$-pair).
Let $S$ and $T$ be proper subsets of $L$. A pair $\left (S,T \right )$ is called a maximal $ST$-pair of $L$ if 
\begin{enumerate}
\item $L$ is $ST$-distributive.
\item There is no proper $S^{'} \supsetneq S$ such that $L$ is $S^{'}T$-distributive.
\item There is no proper $T^{'} \supsetneq T$ such that $L$ is $ST^{'}$-distributive.
\end{enumerate}
\end{definition}

\begin{example} (Maximal and Non-maximal $ST$-pairs). \label{max_stp_ex}
Consider the pair of (proper) subsets $S = \left \{u,v\right \}$ and $T = \left \{w,0\right \}$ from Example \ref{st-distr-ex}. Although, $\mathbf{N}_5$ is $ST$-distributive, this pair  is not a maximal $ST$-pair because if we let $S' =  \left \{0,u,v,w\right \}$ and $T' = \left \{0,1,w\right \}$, then $\mathbf{N}_5$ is $S'T'$-distributive. Now, $\left (S',T' \right )$ is a maximal $ST$-pair of $\mathbf{N}_5$. Clearly, $S'$ cannot be expanded further without making it all $N_5$. On the other hand, trying to add either $u$ or $v$ to $T$ breaks the relative distributive property because $u \land \left( v \lor w \right)$ and $v \lor \left( u \land w \right)$ do not distribute.
\end{example}

We conclude this subsection by remarking that $ST$-distributive lattices are more general than what will be discussed in this paper. At present, we will limit ourselves to finding $ST$-distributive lattices in non-distributive lattices. However, Definition \ref{st-distr} does not require the initial lattice $L$ to be non-distributive. In fact, we mentioned earlier that a lattice $L$ is distributive if and only if it is $ST$-distributive for $S = T = L$. Therefore, the study of $ST$-distributive lattices can encompass problems regarding distributive lattices such as finding and describing their maximal proper sublattices as done in \cite{R}. In such a context, it is helpful to observe that for a proper sublattice $S$ of $L$, the fact that $L$ is an $SS$-distributive lattice is equivalent to $S$ being a distributive (proper) sublattice  of $L$.  Furthermore, if  $\left (S,S \right )$ is a maximal $\left (S,S \right )$-pair, this implies that $S$ is a proper maximal distributive sublattice  of $L$.

\subsection{Basic Properties} \label{basic_props}


With the definitions established, this subsection introduces some basic properties of $ST$-distributive lattices. Two of them, the efficiency criteria and the closure of $S$'s given fixed $T$, will facilitate our search of subsets inducing $ST$-distributivity. \\

\indent The efficiency criteria stated in Property \ref{eff_crit} gives conditions on a triple of elements $a$, $b$, and $c$ of a lattice that each guarantees that the first element under meet (or join) will distributive into the join (or meet) of the other two. The conditions all involve some form of order relation among the elements of the triple. Their advantage is that for some triples, we can determine distribution from simply comparing the elements rather than computing meets and joins.

\begin{property} (Efficiency Criteria). \label{eff_crit}
Let $L$ be a lattice with $a, b, c \in L$ such that the triple $(a,b,c)$ satisfies any of the following conditions:
\begin{center}
\begin{varwidth}{3cm}
\begin{enumerate}
\item $b \leq c$ or $b \geq c$
\item $a \geq b$ and $a \geq c$
\item $a \leq b$ and $a \leq c$ 
\end{enumerate}
\end{varwidth}
\end{center}
Then the following hold:
\begin{eqnarray}
a \land \left(b \lor c\right) = \left(a \land b\right) \lor \left(a \land c\right) \; \:  \label{meet_distr} \\
a \lor \left(b \land c\right) = \left(a \lor b\right) \land \left(a \lor c\right) \; . \label{join_distr}
\end{eqnarray}
\begin{proof}
We show only that each condition implies Equation (\ref{meet_distr}).
Duality gives Equation (\ref{join_distr}) because condition 1 is self-dual and conditions 2 and 3 are each other's duals.\\

\underline{Condition 1}: By commutativity, it is sufficient to show the result for only one of the two inequalities of condition 1. Without loss of generality, suppose that $b \leq c$. Then $a \land b \leq a \land c$. Applying the Connecting Lemma \cite{Bi,DP}, we get
\begin{equation}
a \land \left(b \lor c\right) = a \land c =  \left(a \land b\right) \lor \left(a \land c\right) \; .
\end{equation}

\underline{Condition 2}: Suppose that $a \geq b$ and $a \geq c$. Then by definition of $\lor$ as the supremum, $a \geq b \lor c$. This implies that
\begin{equation} \label{eff_crit_proof_cond2_side1}
a \land \left(b \lor c\right) = b \lor c \; .
\end{equation}
On the other hand, $a \geq b$ and $a \geq c$ also imply that $a \land b = b$ and $a \land c = c$. Therefore,
\begin{equation} \label{eff_crit_proof_cond2_side2}
\left(a \land b\right) \lor \left(a \land c\right) = b \lor c 
\end{equation}
and distribution follows from Equations (\ref{eff_crit_proof_cond2_side1}) and (\ref{eff_crit_proof_cond2_side2}). \\

\underline{Condition 3}: This is similar to the proof for condition 2 except that we get 
\begin{equation}
a \land \left(b \lor c\right) = a = a \lor a = \left(a \land b\right) \lor \left(a \land c\right) \; .
\end{equation}
\qed 
\end{proof}
\end{property}

We make a few comments on Property \ref{eff_crit}. Considering conditions 2 and 3, it is natural to ask what happens if $a$ is comparable to both $b$ and $c$ but in different ways. Then we have either $b \leq a \leq c$ or $c \leq a \leq b$ and hence, condition 1 by transitivity of $\leq$. Thus, it is not necessary to list this as a separate condition. Also note that as a corollary of Property \ref{eff_crit}, we get that we can always add 0 and 1 to any $S$ or $T$ without affecting the relative distributive property. If $a = 0$ or $a = 1$, condition 2 or 3 is satisfied and if either $b$ or $c$ is 0 or 1, then condition 1 is satisfied. In addition, we obtain the following useful fact.

\begin{example} (T-chain). \label{t_chain_gen}
Let $L$ be any lattice (finite or infinite), $S$ any subset of $L$, and $T$ any chain of $L$. Then $L$ is $ST$-distributive by condition 1 of Property \ref{eff_crit}. 
\end{example}

\indent The second property we establish says that if we fix the set $T$ of a pair $\left( S,T \right)$ that we are considering in a lattice $L$, then the set of possible subsets $S$ that make $L$ $ST$-distributive is closed under union. 

\begin{property} (Closure of Union of $S$'s with Fixed $T$). \label{s-closed}
Let $L$ be a lattice with $S_1, S_2, T\subseteq L$. If $L$ is $S_1 T$-distributive and $S_2 T$-distributive, then $L$ is $\left(S_1 \cup S_2\right)T$-distributive.
\begin{proof}
Let $s \in S_1 \cup S_2$. Then $s \in S_1$ or $s \in S_2$. If $s \in S_1$, then for all $t_1, t_2 \in T$,
\begin{eqnarray}
&s \land \left(t_1 \lor t_2\right) &= \left(s \land t_1\right) \lor \left(s \land t_2\right) \label{st-meet} \\
&s \lor \left(t_1 \land t_2\right) &= \left(s \lor t_1\right) \land \left(s \lor t_2\right) \label{st-join}
\end{eqnarray}
because $L$ is $S_1 T$-distributive. The case that $s \in S_2$ is done in the same way. Therefore, $L$ is $\left(S_1 \cup S_2\right)T$-distributive.
\qed
\end{proof}
\end{property}

An important consequence of Property \ref{s-closed} is that for a given $T$, we can form a maximum $S$ such that $L$ is $ST$-distributive by checking all elements in $L$, and picking the ones that distribute into any pair of elements of $T$. This avoids checking each subset of  $L$ and will play a crucial role in computing $ST$-distributivity in the next subsection. \\

\indent We now present an interesting  algebraic extension of the result of Property \ref{s-closed}.

\begin{lemma} (Complete Lattice of $S$'s with Fixed $T$). \label{compl_lat_S}
Let $L$ be any lattice. For a fixed $T$, the collection of all subsets $S \subseteq L$ such that $L$ is $ST$-distributive forms a complete distributive lattice of sets. The join and meet operations are given by union and intersection respectively. We denote this lattice $\mathcal{S} \left (L,T \right )$.
\begin{proof}
Similar to Property \ref{s-closed}. \qed
\end{proof}
\end{lemma}

When $L$ is a finite lattice, the lattice $\mathcal{S} \left (L,T \right )$ is also finite. In this case, the latter is a bounded lattice whose top element is the set of all elements in $L$ that are meet and join distributive in $T$. This is
$L$ itself if $L$ is distributive. When $L$ is not distributive, it is a subset of $L$ the contains the desired set of step 1 of Algorithm \ref{st-algo} (it may be larger since the algorithm has additional restrictions from Problem \ref{st-distr-search}). This subset may or may not be proper (e.g. it is all $L$ if $T$ is a chain).\\

\indent One last matter that we touch is how to construct new $ST$-distributive lattices from known ones. It is known that the distributive property of lattices is preserved by different lattice constructs: sublattices, products of lattices, lattice homomorphisms, and dual lattices. The following proposition 
generalizes this to $ST$-distributive lattices. The proofs come straight from the definitions and are thus omitted.

\begin{proposition} (Preservation of $ST$-distributivity by Lattice Constructs).
Let $L$ be a lattice with subsets $S$ and $T$ such that $L$ is $ST$-distributive.
\begin{enumerate}
\item If $K$ is a sublattice of $L$, then $K$ is $S_1 T_1$-distributive for $S_1 = S \cap K$ and $T_1 = T \cap K$.
\item If $L'$ is another lattice with subsets $S'$ and $T'$ such that $L'$ is $S'T'$-distributive, then $L \times L'$ is $\left (S \times S' \right ) \left (T \times T' \right )$-distributive.
\item If $\phi: L \rightarrow K$ is a lattice homomorphism that is onto, then $K$ is $\phi \left (S \right ) \phi \left (T \right )$-distributive.
\item If $L^{\partial}$ represents the dual lattice of $L$ and $S$ and $T$ are sublattices with duals $S^{\partial}$ and $T^{\partial}$, then $L^{\partial}$ is $S^{\partial}T^{\partial}$-distributive.
\end{enumerate}
\end{proposition}

\subsection{Search for $ST$-distributive Lattices} \label{search}


We want to study non-distributive lattices to find the subsets $S$ and $T$ for which they are $ST$-distributive. The long-term goal is to identify families of lattices with a particular characterization of the subsets inducing this. As a result, we have the following initial problem, where $L^{*}$ denotes the set of non-identity elements of a lattice $L$, i.e. $L^{*} = L \setminus \left \{0,1\right \}$.

\begin{problem} ($ST$-distributive Search). \label{st-distr-search}
Given a non-distributive lattice $L$, we want to find all of its maximal $ST$-pairs subject to the following conditions:
\begin{center}
\begin{varwidth}{3.5cm}
\begin{enumerate}
\item $S \cap T = \emptyset$
\item $S, T \subseteq L^{*}$
\item $S \neq \emptyset$ and $T \neq \emptyset$
\end{enumerate}
\end{varwidth}
\end{center}
\end{problem}

We explain the reasons for our conditions. The disjointness of $S$ and $T$ is for narrowing the search space. Our disregard for 0 and 1 stems from the fact that they can always be added to any $S$ or $T$ without affecting the relative distributive property of the pair. This is a consequence of Property \ref{eff_crit} mentioned earlier. Finally, the non-emptiness of both $S$ and $T$ is to focus only on non-trivial pairs with actual distribution of elements.\\

We make a \emph{very important clarification} before moving forward. \emph{From now on, when we refer to a maximal $ST$-pair of a lattice $L$, we will mean a pair $\left (S, T \right )$ that satisfies not only Definition \ref{max_stp} but also the conditions of Problem \ref{st-distr-search}}. In particular, we remark that the maximality of a set of the pair will be altered by both the disjointness condition and the disregard of proper subsets of $L$ not contained in $L^{*}$. \\

\indent Having established the specifications of our search, we explain our method for carrying in out. We write a SageMath program that does an ``intelligent" exhaustive search by applying Properties \ref{eff_crit} and \ref{s-closed}. An overview of how it works is given in Algorithm \ref{st-algo}. Property \ref{s-closed} allows us to construct a maximum $S$ for each $T$ by iterating across all elements of $L^{*} \setminus T$, putting those that distribute in one set. The resulting pair $\left (S,T \right )$ is a maximal $ST$-pair unless there is a $T'$ strictly containing $T$ that has the same maximum $S$. This motivates our processing of the subsets $T$ by increasing size. Property \ref{eff_crit} reduces the computations done when determining who can be added to $S$. After running the program, we manually verify the results for correctness and completeness. This finishes our discussion of our search problem.

\begin{algorithm} ($ST$-program Algorithm). \label{st-algo}
\begin{itemize}
\item \underline{Specifications}: \\
\begin{enumerate}
\item Input: lattice $L$ 
\item Output: list of maximal $ST$-pairs of $L$ with Problem \ref{st-distr-search} conditions\\
\end{enumerate} 

\item \underline{Process Overview}:\\

For each $T \subseteq L^{*}$ (by increasing size):
\begin{enumerate}
\item Build largest possible set $S$ such that $L$ is $ST$-distributive. 
\begin{enumerate}
\item Add to candidate set $S$ all elements $s$ of $L^{*} \setminus T$ such that $s$ distributes into $T$.
\end{enumerate}
\item If returned $S$ is non-empty, add pair $\left(S,T\right)$ to list of pairs and remove any pairs contained by $\left(S,T\right)$.
\end{enumerate}
\end{itemize}
\end{algorithm}

\subsection{Family $\mathbf{M}_{n,n}$} \label{fam-M_n,n}


Before we begin, we clarify some notation. For any natural number $n$, $\mathbf{n}$ denotes the chain of $n$ elements and $\mathbf{\bar{n}}$ the antichain of $n$ elements, respectively. In addition, for lattices $L$ and $K$, $L \oplus K$ denotes their linear sum and $L \times K$ their product.\\

\indent We now introduce the family of lattices with which we have begun our search of $ST$-distributive lattices: $\mathbf{M}_{n,n}$.  It consists of lattices obtained by gluing two copies of $\mathbf{M}_n = \mathbf{1} \oplus \mathbf{\bar{n}} \oplus \mathbf{1}$ in a particular manner for each natural number $n \geq 3$. The pictorial idea is to take the rightmost edge from the top element in the Hasse diagram of one of the copies and make it equal to the leftmost edge from the bottom element of the other copy's diagram. Definition \ref{M_n,n} formally describes the family and presents relevant notation. Some members of the family are given in Example \ref{M_n,n-ex}.

\begin{definition} (Lattice $\mathbf{M}_{n,n}$). \label{M_n,n}
For a natural number $n \geq 3$, the lattice $\mathbf{M}_{n,n}$ has element set $\left \{0,1,a_1, ..., a_n, b_1, ..., b_n\right \}$ with $0$ and $1$ being the bottom and top element respectively. $0$ is covered by all of the $a_i$'s while $1$ covers all of the $b_j$'s. In addition $a_n$ is covered by all the $b_j$'s and $b_1$ covers all of the $a_i$'s. No other covering relations exist. This results in $\mathbf{M}_{n,n}$ having two isomorphic copies of $\mathbf{M}_n$: $\left \{0, a_1, a_2, ..., a_n, b_1\right \}$ and  $\left \{a_n, b_1, b_2, ..., b_n, 1 \right \}$. Note that both of them share $a_n$ and $b_1$ (the shared edge from the pictorial description above).
\end{definition}

\begin{example} (Examples of $\mathbf{M}_{n,n}$). \label{M_n,n-ex}
$\mathbf{M}_{3,3}$ and $\mathbf{M}_{4,4}$ are shown in Figs. \ref{m33} and \ref{m44} respectively. Note that 

\begin{eqnarray}
\mathbf{M}_3 \cong \left \{0, a_1, a_2, a_3, b_1\right \} \cong \left \{a_3, b_1, b_2, b_3, 1 \right \} 
\quad \; \; \,\\
\mathbf{M}_4 \cong \left \{0, a_1, a_2, a_3, a_4, b_1\right \} \cong \left \{a_4, b_1, b_2, b_3, b_4, 1 \right \}.
\end{eqnarray}

\end{example}

 \begin{figure}
 \begin{minipage}{0.45\textwidth}
 \centering
 \begin{tikzpicture}
[roundnode/.style={circle, draw=black, fill=white, very thin, text width =4.00mm, inner sep = 0pt, text centered}]
 \node [roundnode] (1) at (.5,1.5) {$1$};
  \node [roundnode] (a) at (-.5,0.5) {$b_1$};
  \node [roundnode] (b) at (.5,0.5) {$b_2$};
  \node [roundnode] (c) at (1.5,0.5) {$b_3$};
   \node [roundnode] (d) at (-1.5,-0.5) {$a_1$};
  \node [roundnode] (e) at (-.5,-0.5) {$a_2$};
  \node [roundnode] (f) at (.5,-0.5) {$a_3$};
  \node [roundnode] (0) at (-.5,-1.5) {$0$};
 \draw (1) -- (a) -- (d) -- (0) -- (f) -- (c) -- (1);
 \draw (1) -- (b) -- (f) -- (a) -- (e) -- (0);
\end{tikzpicture}
\caption{Lattice $\mathbf{M}_{3,3}$} \label{m33}
 \end{minipage}
 \begin{minipage}{0.45\textwidth}
 \centering
 \begin{tikzpicture}
[roundnode/.style={circle, draw=black, fill=white, very thin, text width =4.00mm, inner sep = 0pt, text centered}]
 \node [roundnode] (1) at (.5,1.5) {$1$};
  \node [roundnode] (b1) at (-.5,0.5) {$b_1$};
  \node [roundnode] (b2) at (0.15,0.5) {$b_2$};
  \node [roundnode] (b3) at (0.85,0.5) {$b_3$};
   \node [roundnode] (b4) at (1.5,0.5) {$b_4$};
   \node [roundnode] (a1) at (-1.5,-0.5) {$a_1$};
  \node [roundnode] (a2) at (-.85,-0.5) {$a_2$};
   \node [roundnode] (a3) at (-.15,-0.5) {$a_3$};
  \node [roundnode] (a4) at (.5,-0.5) {$a_4$};
  \node [roundnode] (0) at (-.5,-1.5) {$0$};
 \draw (1) -- (b1) -- (a1) -- (0) -- (a2) -- (b1) -- (a3) -- (0) -- (a4) -- (b1);
 \draw (1) -- (b2) -- (a4) -- (b3) -- (1) -- (b4) -- (a4);
\end{tikzpicture}
\caption{Lattice $\mathbf{M}_{4,4}$} \label{m44}
 \end{minipage}
\end{figure}


We briefly discuss the family $\mathbf{M}_{n,n}$. First, we introduce some terms we use to refer to particular subsets and elements of $\mathbf{M}_{n,n}$. This vocabulary will facilitate the discussion of our results in the next subsection.

\begin{definition} (Level).  \label{level}
In $\mathbf{M}_{n,n}$, a level is a set of elements that constitutes the antichain ($\mathbf{\bar{n}}$) of a copy of $\mathbf{M}_n$. In other words, the sets $\left \{a_1, ..., a_n\right \}$ and $\left \{b_1, ..., b_n\right \}$.
\end{definition}

\begin{definition} (Link Element). \label{link}
In $\mathbf{M}_{n,n}$, a link element (or link) is one of the two elements shared by both isomorphic copies of $\mathbf{M}_n$, that is, $a_n$ and $b_1$.
\end{definition}


We complete our introduction to the family $\mathbf{M}_{n,n}$ by remarking that the lattice $\mathbf{M}_{n,n}$ is both modular and non-distributive for $n \geq 3$. Both facts follow from applying the $\mathbf{M}_3$-$\mathbf{N}_5$ Theorem \cite{Bi,DP,G2011}: $\mathbf{M}_{n,n}$ always has two isomorphic copies of $\mathbf{M}_n$ as previously mentioned (and hence of $\mathbf{M}_3$ since $n \geq 3$),  but never an isomorphic copy of $\mathbf{N}_5$. Thus, we have a family of modular and non-distributive lattices. This is good because it allows us to study $ST$-distributivity non-trivially in a relatively ``controlled" environment. Finally, note that it is possible to construct $\mathbf{M}_{2,2}$ and $\mathbf{M}_{1,1}$ following Definition \ref{M_n,n}. However, the resulting lattices would be isomorphic to $\mathbf{2} \times \mathbf{3}$ and $\mathbf{4}$ ($\mathbf{n}$ denotes the chain of $n$ elements) and thus, distributive. This is what motivates our specification that $n \geq 3$ when defining the family.

\subsection{Results} \label{results}

Here, we present our results in the study of $ST$-distributivity in the family $\mathbf{M}_{n,n}$. We start with some computational results and then present theoretical results derived from them, concluding with a complete characterization of the maximal $ST$-pairs of this lattice family. Throughout the whole discussion, the reader should bear in mind the clarification we made in Subsect. \ref{search} about what we mean when we refer to a maximal $ST$-pair.

\subsubsection{Computational Results.}
\indent Running the program discussed in Algorithm \ref{st-algo} with inputs $\mathbf{M}_{3,3}$ and $\mathbf{M}_{4,4}$ yields 27 and 42 maximal $ST$-pairs respectively. These are listed in Tables \ref{m33-pairs} and \ref{m44-pairs}. Manual verification then shows that these lists are correct and complete. Further study of these pairs results in identifying 5 types into which all of them belong. These types are described with constructions that select elements based on their role in the structure of  $\mathbf{M}_{n,n}$ (links and levels).

\begin{remark} (5 Types of Maximal $ST$-pairs). \label{pair-types}
All of the maximal $ST$-pairs of $\mathbf{M}_{3,3}$ and $\mathbf{M}_{4,4}$ are of one of the following 5 types. The type of each pair indicated in Tables \ref{m33-pairs} and \ref{m44-pairs}. An example of each type in $\mathbf{M}_{4,4}$ is illustrated in Figures \ref{m44_t-chain}--\ref{m44_s-level-minus-link} ($S$ in cyan, $T$ in black).
\begin{enumerate}
\item \underline{\emph{T-chain}}: $T$ is a chain. $S$ has all other elements. 
\item \underline{\emph{S-link}}: $S$ has only one of the two links. $T$ has one more element in this link's level and all of the elements in the other level.
\item \underline{\emph{S-2-links}}: $S$ has both links. $T$ has one non-link element of each level
\item \underline{\emph{S-level}}: $S$ is one of the two levels. $T$ has two elements in the other level, one of which must be its link.
\item \underline{\emph{S-level-minus-link}}: $S$ is one of the two levels minus its link. $T$ is both links plus another element from the level that does not contain $S$.
\end{enumerate}

\end{remark}


\begin{table}
\caption{Maximal $ST$-pairs of $\mathbf{M}_{3,3}$} \label{m33-pairs}
\begin{center}
\begin{tabular}{ @{\makebox[3em][c]{\rownumber\space}}  c c c }
\hline
   $S$ & $T$ & Type
  \gdef\rownumber{\stepcounter{magicrownumbers}\arabic{magicrownumbers}} \\
  \hline
    $\left \{a_2,a_3,b_1,b_2,b_3\right \}$ & $\left \{a_1\right \}$ & 1\\
  \hline
  $\left \{a_1,a_3,b_1,b_2,b_3\right \}$ & $\left \{a_2\right \}$ & 1\\
  \hline
   $\left \{a_1,a_2,b_1,b_2,b_3\right \}$ & $\left \{a_3\right \}$ & 1\\
  \hline
   $\left \{a_1,a_2,a_3,b_2,b_3\right \}$ & $\left \{b_1\right \}$ & 1\\
  \hline
  $\left \{a_1,a_2,a_3,b_1,b_3\right \}$ & $\left \{b_2\right \}$ & 1\\
  \hline
  $\left \{a_1,a_2,a_3,b_1,b_2\right \}$ & $\left \{b_3\right \}$ & 1\\
  \hline
  $\left \{a_2,a_3,b_2,b_3\right \}$ & $\left \{a_1,b_1\right \}$ & 1\\
  \hline
  $\left \{a_1,a_3,b_2,b_3\right \}$ & $\left \{a_2,b_1\right \}$ & 1\\
  \hline
  $\left \{a_1,a_2,b_2,b_3\right \}$ & $\left \{a_3,b_1\right \}$ & 1\\
  \hline
  $\left \{a_1,a_2,b_1,b_3\right \}$ & $\left \{a_3,b_2\right \}$ & 1\\
  \hline
  $\left \{a_1,a_2,b_1,b_2\right \}$ & $\left \{a_3,b_3\right \}$ & 1\\
   \hline
  $\left \{b_1\right \}$ & $\left \{a_1,a_2,a_3,b_2\right \}$ & 2\\
  \hline
  $\left \{b_1\right \}$ & $\left \{a_1,a_2,a_3,b_3\right \}$ & 2\\
  \hline
   $\left \{a_3\right \}$ & $\left \{a_1,b_1,b_2,b_3\right \}$ & 2\\
  \hline
  $\left \{a_3\right \}$ & $\left \{a_2,b_1,b_2,b_3\right \}$ & 2\\
  \hline
 $\left \{a_3,b_1\right \}$ & $\left \{a_1,b_2\right \}$ & 3\\
  \hline
  $\left \{a_3,b_1\right \}$ & $\left \{a_1,b_3\right \}$ & 3\\
   \hline
  $\left \{a_3,b_1\right \}$ & $\left \{a_2,b_2\right \}$ & 3\\
   \hline
  $\left \{a_3,b_1\right \}$ & $\left \{a_2,b_3\right \}$ & 3\\
   \hline
  $\left \{a_1,a_2, a_3\right \}$ & $\left \{b_1,b_2\right \}$ & 4\\
  \hline
  $\left \{a_1,a_2, a_3\right \}$ & $\left \{b_1,b_3\right \}$ & 4\\
  \hline
  $\left \{b_1,b_2,b_3\right \}$ & $\left \{a_1,a_3\right \}$ & 4\\
  \hline
  $\left \{b_1,b_2,b_3\right \}$ & $\left \{a_2,a_3\right \}$ & 4\\
  \hline
   $\left \{a_1,a_2\right \}$ & $\left \{a_3,b_1,b_2\right \}$ & 5\\
   \hline
 $\left \{a_1,a_2\right \}$ & $\left \{a_3,b_1,b_3\right \}$ & 5\\
   \hline
  $\left \{b_2,b_3\right \}$ & $\left \{a_1,a_3,b_1\right \}$ & 5\\
  \hline
  $\left \{b_2,b_3\right \}$ & $\left \{a_2,a_3,b_1\right \}$ & 5\\
  \hline
\end{tabular}
\end{center}
\end{table}

\begin{table}
\caption{Maximal $ST$-pairs of $\mathbf{M}_{4,4}$} \label{m44-pairs}
\begin{center}
\begin{tabular}{ @{\makebox[3em][c]{\foo\space}}  c  c  c}
\hline
   $S$ & $T$ & Type 
   \gdef\foo{\stepcounter{magicrownumbers2}\arabic{magicrownumbers2}}  \\
  \hline
  $\left \{a_2,a_3,a_4,b_1,b_2,b_3,b_4\right \}$ & $\left \{a_1\right \}$ & 1\\
    \hline
  $\left \{a_1,a_3,a_4,b_1,b_2,b_3,b_4\right \}$ & $\left \{a_2\right \}$ & 1\\
   \hline
  $\left \{a_1,a_2,a_4,b_1,b_2,b_3,b_4\right \}$ & $\left \{a_3\right \}$ & 1\\
    \hline
  $\left \{a_1,a_2,a_3,b_1,b_2,b_3,b_4\right \}$ & $\left \{a_4\right \}$ & 1\\
    \hline
  $\left \{a_1,a_2,a_3,a_4,b_2,b_3,b_4\right \}$ & $\left \{b_1\right \}$ & 1\\
   \hline
  $\left \{a_1,a_2,a_3,a_4,b_1,b_3,b_4\right \}$ & $\left \{b_2\right \}$ & 1\\
     \hline
  $\left \{a_1,a_2,a_3,a_4,b_1,b_2,b_4\right \}$ & $\left \{b_3\right \}$ & 1\\
     \hline
  $\left \{a_1,a_2,a_3,a_4,b_1,b_2,b_3\right \}$ & $\left \{b_4\right \}$ & 1\\
  
   \hline
  $\left \{a_2,a_3,a_4,b_2,b_3,b_4\right \}$ & $\left \{a_1,b_1\right \}$ & 1\\
  \hline
  $\left \{a_1,a_3,a_4,b_2,b_3,b_4\right \}$ & $\left \{a_2,b_1\right \}$ & 1\\
  \hline
  $\left \{a_1,a_2,a_4,b_2,b_3,b_4\right \}$ & $\left \{a_3,b_1\right \}$ & 1\\
      \hline
  $\left \{a_1,a_2,a_3,b_2,b_3,b_4\right \}$ & $\left \{a_4,b_1\right \}$ & 1\\  
     \hline
  $\left \{a_1,a_2,a_3,b_1,b_3,b_4\right \}$ & $\left \{a_4,b_2\right \}$ & 1\\
     \hline
  $\left \{a_1,a_2,a_3,b_1,b_2,b_4\right \}$ & $\left \{a_4,b_3\right \}$ & 1\\
  \hline
  $\left \{a_1,a_2,a_3,b_1,b_2,b_3\right \}$ & $\left \{a_4,b_4\right \}$ & 1\\
  
   \hline
  $\left \{b_1\right \}$ & $\left \{a_1,a_2,a_3,a_4,b_2\right \}$ & 2\\
   \hline
  $\left \{b_1\right \}$ & $\left \{a_1,a_2,a_3,a_4,b_3\right \}$ & 2\\
   \hline
  $\left \{b_1\right \}$ & $\left \{a_1,a_2,a_3,a_4,b_4\right \}$ & 2\\
   \hline
  $\left \{a_4\right \}$ & $\left \{a_1,b_1,b_2,b_3,b_4\right \}$ & 2\\
  \hline
  $\left \{a_4\right \}$ & $\left \{a_2,b_1,b_2,b_3,b_4\right \}$ & 2\\
  \hline
  $\left \{a_4\right \}$ & $\left \{a_3,b_1,b_2,b_3,b_4\right \}$ & 2\\

   \hline
  $\left \{a_4,b_1\right \}$ & $\left \{a_1,b_2\right \}$ & 3\\
   \hline
  $\left \{a_4,b_1\right \}$ & $\left \{a_1,b_3\right \}$ & 3\\
   \hline
  $\left \{a_4,b_1\right \}$ & $\left \{a_1,b_4\right \}$ & 3\\
    \hline
  $\left \{a_4,b_1\right \}$ & $\left \{a_2,b_2\right \}$ & 3\\
   \hline
  $\left \{a_4,b_1\right \}$ & $\left \{a_2,b_3\right \}$ & 3\\
   \hline
  $\left \{a_4,b_1\right \}$ & $\left \{a_2,b_4\right \}$ & 3\\  
    \hline
  $\left \{a_4,b_1\right \}$ & $\left \{a_3,b_2\right \}$ & 3\\
    \hline
  $\left \{a_4,b_1\right \}$ & $\left \{a_3,b_3\right \}$ & 3\\
   \hline
  $\left \{a_4,b_1\right \}$ & $\left \{a_3,b_4\right \}$ & 3\\
  
   \hline
  $\left \{a_1,a_2,a_3,a_4\right \}$ & $\left \{b_1,b_2\right \}$ & 4\\
   \hline
  $\left \{a_1,a_2,a_3,a_4\right \}$ & $\left \{b_1,b_3\right \}$ & 4\\
     \hline
  $\left \{a_1,a_2,a_3,a_4\right \}$ & $\left \{b_1,b_4\right \}$ & 4\\
    \hline
  $\left \{b_1,b_2,b_3,b_4\right \}$ & $\left \{a_1,a_4\right \}$ & 4\\
    \hline
  $\left \{b_1,b_2,b_3,b_4\right \}$ & $\left \{a_2,a_4\right \}$ & 4\\
    \hline
  $\left \{b_1,b_2,b_3,b_4\right \}$ & $\left \{a_3,a_4\right \}$ & 4\\
  
   \hline
  $\left \{a_1,a_2,a_3\right \}$ & $\left \{a_4,b_1,b_2\right \}$ & 5\\
   \hline
  $\left \{a_1,a_2,a_3\right \}$ & $\left \{a_4,b_1,b_3\right \}$ & 5\\
   \hline
  $\left \{a_1,a_2,a_3\right \}$ & $\left \{a_4,b_1,b_4\right \}$ & 5\\
    \hline
  $\left \{b_2,b_3,b_4\right \}$ & $\left \{a_1,a_4,b_1\right \}$ & 5\\
   \hline
  $\left \{b_2,b_3,b_4\right \}$ & $\left \{a_2,a_4,b_1\right \}$ & 5\\
   \hline
  $\left \{b_2,b_3,b_4\right \}$ & $\left \{a_3,a_4,b_1\right \}$ & 5\\
  \hline
\end{tabular}
\end{center}
\end{table}


 \begin{figure}
 \begin{minipage}{0.32\textwidth}
 \centering
 \begin{tikzpicture}
[roundnode/.style={circle, draw=black, fill=white, very thin, text width =4.00mm, inner sep = 0pt, text centered}, 
snode/.style={circle, draw=black, fill=cyan, very thin, text width =4.00mm, inner sep = 0pt, text centered},
tnode/.style={circle, draw=black, fill=black, very thin, text width =4.00mm, inner sep = 0pt, text centered, text = white}]
 \node [roundnode] (1) at (.5,1.5) {$1$};
 \node [tnode] (b1) at (-.5,0.5) {$b_1$};
  \node [snode] (b2) at (0.15,0.5) {$b_2$};
  \node [snode] (b3) at (0.85,0.5) {$b_3$};
   \node [snode] (b4) at (1.5,0.5) {$b_4$};
   \node [tnode] (a1) at (-1.5,-0.5) {$a_1$};
  \node [snode] (a2) at (-.85,-0.5) {$a_2$};
   \node [snode] (a3) at (-.15,-0.5) {$a_3$};
  \node [snode] (a4) at (.5,-0.5) {$a_4$};
  \node [roundnode] (0) at (-.5,-1.5) {$0$};
 \draw (1) -- (b1) -- (a1) -- (0) -- (a2) -- (b1) -- (a3) -- (0) -- (a4) -- (b1);
 \draw (1) -- (b2) -- (a4) -- (b3) -- (1) -- (b4) -- (a4);
\end{tikzpicture}
\caption{T-chain pair} \label{m44_t-chain}
 \end{minipage}
 \begin{minipage}{0.32\textwidth}
 \centering
\begin{tikzpicture}
[roundnode/.style={circle, draw=black, fill=white, very thin, text width =4.00mm, inner sep = 0pt, text centered}, 
snode/.style={circle, draw=black, fill=cyan, very thin, text width =4.00mm, inner sep = 0pt, text centered},
tnode/.style={circle, draw=black, fill=black, very thin, text width =4.00mm, inner sep = 0pt, text centered, text = white}]
 \node [roundnode] (1) at (.5,1.5) {$1$};
 \node [snode] (b1) at (-.5,0.5) {$b_1$};
 \node [tnode] (b2) at (0.15,0.5) {$b_2$};
 \node [roundnode] (b3) at (0.85,0.5) {$b_3$};
 \node [roundnode] (b4) at (1.5,0.5) {$b_4$};
  \node [tnode] (a1) at (-1.5,-0.5) {$a_1$};
\node [tnode] (a2) at (-.85,-0.5) {$a_2$};
 \node [tnode] (a3) at (-.15,-0.5) {$a_3$};
\node [tnode] (a4) at (.5,-0.5) {$a_4$};
\node [roundnode] (0) at (-.5,-1.5) {$0$};
\draw (1) -- (b1) -- (a1) -- (0) -- (a2) -- (b1) -- (a3) -- (0) -- (a4) -- (b1);
\draw (1) -- (b2) -- (a4) -- (b3) -- (1) -- (b4) -- (a4);
\end{tikzpicture}
\caption{S-link pair} \label{m44_s-link}
 \end{minipage}
 \begin{minipage}{0.32\textwidth}
 \centering
 \begin{tikzpicture}
[roundnode/.style={circle, draw=black, fill=white, very thin, text width =4.00mm, inner sep = 0pt, text centered}, 
snode/.style={circle, draw=black, fill=cyan, very thin, text width =4.00mm, inner sep = 0pt, text centered},
tnode/.style={circle, draw=black, fill=black, very thin, text width =4.00mm, inner sep = 0pt, text centered, text = white}]
 \node [roundnode] (1) at (.5,1.5) {$1$};
 \node [snode] (b1) at (-.5,0.5) {$b_1$};
  \node [tnode] (b2) at (0.15,0.5) {$b_2$};
  \node [roundnode] (b3) at (0.85,0.5) {$b_3$};
   \node [roundnode] (b4) at (1.5,0.5) {$b_4$};
   \node [tnode] (a1) at (-1.5,-0.5) {$a_1$};
  \node [roundnode] (a2) at (-.85,-0.5) {$a_2$};
   \node [roundnode] (a3) at (-.15,-0.5) {$a_3$};
  \node [snode] (a4) at (.5,-0.5) {$a_4$};
  \node [roundnode] (0) at (-.5,-1.5) {$0$};
 \draw (1) -- (b1) -- (a1) -- (0) -- (a2) -- (b1) -- (a3) -- (0) -- (a4) -- (b1);
 \draw (1) -- (b2) -- (a4) -- (b3) -- (1) -- (b4) -- (a4);
\end{tikzpicture}
\caption{S-2-links pair} \label{m44_s-2-links}
 \end{minipage}

\vspace{3mm}


 \begin{minipage}{0.45\textwidth}
 \centering
\begin{tikzpicture}
[roundnode/.style={circle, draw=black, fill=white, very thin, text width =4.00mm, inner sep = 0pt, text centered}, 
snode/.style={circle, draw=black, fill=cyan, very thin, text width =4.00mm, inner sep = 0pt, text centered},
tnode/.style={circle, draw=black, fill=black, very thin, text width =4.00mm, inner sep = 0pt, text centered, text = white}]
 \node [roundnode] (1) at (.5,1.5) {$1$};
 \node [tnode] (b1) at (-.5,0.5) {$b_1$};
 \node [tnode] (b2) at (0.15,0.5) {$b_2$};
 \node [roundnode] (b3) at (0.85,0.5) {$b_3$};
 \node [roundnode] (b4) at (1.5,0.5) {$b_4$};
  \node [snode] (a1) at (-1.5,-0.5) {$a_1$};
\node [snode] (a2) at (-.85,-0.5) {$a_2$};
 \node [snode] (a3) at (-.15,-0.5) {$a_3$};
\node [snode] (a4) at (.5,-0.5) {$a_4$};
\node [roundnode] (0) at (-.5,-1.5) {$0$};
\draw (1) -- (b1) -- (a1) -- (0) -- (a2) -- (b1) -- (a3) -- (0) -- (a4) -- (b1);
\draw (1) -- (b2) -- (a4) -- (b3) -- (1) -- (b4) -- (a4);
\end{tikzpicture}
\caption{S-level pair} \label{m44_s-level}
 \end{minipage}
  \begin{minipage}{0.45\textwidth}
 \centering
 \begin{tikzpicture}
[roundnode/.style={circle, draw=black, fill=white, very thin, text width =4.00mm, inner sep = 0pt, text centered}, 
snode/.style={circle, draw=black, fill=cyan, very thin, text width =4.00mm, inner sep = 0pt, text centered},
tnode/.style={circle, draw=black, fill=black, very thin, text width =4.00mm, inner sep = 0pt, text centered, text = white}]
 \node [roundnode] (1) at (.5,1.5) {$1$};
 \node [tnode] (b1) at (-.5,0.5) {$b_1$};
 \node [tnode] (b2) at (0.15,0.5) {$b_2$};
 \node [roundnode] (b3) at (0.85,0.5) {$b_3$};
 \node [roundnode] (b4) at (1.5,0.5) {$b_4$};
  \node [snode] (a1) at (-1.5,-0.5) {$a_1$};
\node [snode] (a2) at (-.85,-0.5) {$a_2$};
 \node [snode] (a3) at (-.15,-0.5) {$a_3$};
\node [tnode] (a4) at (.5,-0.5) {$a_4$};
\node [roundnode] (0) at (-.5,-1.5) {$0$};
\draw (1) -- (b1) -- (a1) -- (0) -- (a2) -- (b1) -- (a3) -- (0) -- (a4) -- (b1);
\draw (1) -- (b2) -- (a4) -- (b3) -- (1) -- (b4) -- (a4);
\end{tikzpicture}
\caption{S-level-minus-link pair} \label{m44_s-level-minus-link}
 \end{minipage}
\end{figure}

The natural questions arising from Remark \ref{pair-types} are if these constructions of maximal $ST$-pairs work for all other $\mathbf{M}_{n,n}$  (recalling $n \geq 3$) and whether larger $n$'s yield maximal $ST$-pairs of other forms. We already know that type T-chain will always work by Example \ref{t_chain_gen}. We show that the remaining constructions work for all lattices of the $\mathbf{M}_{n,n}$ family and prove that these five constructions are the only ones possible in it. 

\subsubsection{Preliminary Theoretical Results.}
\indent To develop our main theoretical results, we require some preliminary results: the modular comparability criteria for distributivity (MCCD) and Proposition \ref{and-m_nn}. The MCCD gives a simple condition that guarantees distribution for any triple of elements in a non-distributive but modular lattice. 

\begin{proposition} (Modular Comparability Criteria for Distributivity --- MCCD). \label{mccd}
Let $L$ be a modular lattice with elements $a,b,c \in L$. If there is any comparability between any 2 of the 3 elements, then
\begin{eqnarray}
a \land \left(b \lor c\right) &= \left(a \land b\right) \lor \left(a \land c\right) \label{meet-distr}  \\
a \lor \left(b \land c\right) &= \left(a \lor b\right) \land \left(a \lor c\right) \label{join-distr}  .
\end{eqnarray} 
\begin{proof}
By duality, we need only to show Eq. (\ref{meet-distr}). There are three cases. First, If $b \leq c$ or $b \geq c$, we are done by Property \ref{eff_crit}. Second, if $a \leq b$ or $a \leq c$, then $a \leq b \lor c$. Thus, the Connecting Lemma \cite{Bi,DP} and the absorption property of lattices give
\begin{equation}
a \land \left(b \lor c\right) = a = \left(a \land b\right) \lor \left(a \land c\right) \; .
\end{equation}
Finally, we consider $a \geq b$ or $a \geq c$. Suppose that $a \geq c$, then by modularity and the Connecting Lemma \cite{Bi,DP},
\begin{equation}
a \land \left(b \lor c\right) = \left(a \land b\right) \lor c = \left(a \land b\right) \lor \left(a \land c\right) \; .
\end{equation}
The sub-case $a \geq b$ then follows from commutativity of $\lor$. Therefore, Equation (\ref{meet-distr}) is satisfied whenever $\left \{a,b,c\right \}$ is not an antichain.
\qed
\end{proof}
\end{proposition}

The importance of the MCCD is that it implies that we need only to worry about triples that are antichains when checking for relative distributivity. Note that in general, the converse of MCCD is not true: a triple that is an antichain may distribute (e.g. a triple of singleton sets in a power set lattice). However,  in $\mathbf{M}_{n,n}$, antichains never distribute.

\begin{proposition} (Antichains Not Distributive in $\mathbf{M}_{n,n}$). \label{and-m_nn}
Suppose $n \geq 3$ and $x,y,z \in \mathbf{M}_{n,n}$ are distinct elements. If $\left \{x,y,z\right \}$ is an antichain, then 
\begin{eqnarray}
x \land \left(y \lor z\right) &\neq \left(x \land y\right) \lor \left(x \land z\right) \\
y \land \left(x \lor z\right) &\neq \left(y \land x\right) \lor \left(y \land z\right) \\
z \land \left(x \lor y\right) &\neq \left(z \land x\right) \lor \left(z \land y\right)  .
\end{eqnarray}
\begin{proof}
There are two ways of choosing $x,y,z \in \mathbf{M}_{n,n}$ such that $\left \{x,y,z\right \}$ is an antichain. \\

\underline{Case 1}: $x,y,z$ are all elements of the same level of $\mathbf{M}_{n,n}$. Without loss of generality, suppose that $\left \{x,y,z\right \} = \left \{a_i, a_j, a_k\right \}$. Then we are done by the fact that $a_i, a_j, a_k$ are three elements in the $\bf{\bar{n}}$ of an isomorphic copy of $\mathbf{M}_n$ in $\mathbf{M}_{n,n}$.\\

\underline{Case 2}: $x,y,z$ consist of three non-link elements split among the two levels of $\mathbf{M}_{n,n}$. Without loss of generality, suppose that $\left \{x,y,z\right \} = \left \{a_i, a_j, b_k\right \}$ with $i,j \neq n$ and $k \neq 1$. Then
\begin{eqnarray}
a_i \land \left(a_j \lor b_k\right) = a_i \land 1 = a_i \neq 0 = 0 \lor 0 = \left(a_i \land a_j\right) \lor \left(a_i \land b_k\right)  \\
a_j \land \left(a_i \lor b_k\right) = a_j \land 1 = a_j \neq 0 = 0 \lor 0 = \left(a_j \land a_i\right) \lor \left(a_j \land b_k\right)  \\
b_k \land \left(a_i \lor a_j\right) = b_k \land b_1 = a_n \neq 0 = 0 \lor 0 = \left(b_k \land a_i\right) \lor \left(b_k \land a_j\right) .
\end{eqnarray}
\qed
\end{proof}
\end{proposition}

The implication of Proposition \ref{and-m_nn} is that for a pair of subsets $S$ and $T$ of a lattice $L$ to be a maximal $ST$-pair, there can be no antichain of three distinct elements in $L$ that has one element in $S$ and two in $T$. This motivates the following definition.

\begin{definition} ($ST$-breaking Antichain).
Given $S, T \subseteq \mathbf{M}_{n,n}^{*}$, an $ST$-breaking antichain is an ordered triple $(x,y,z)$ of elements of $\mathbf{M}_{n,n}^{*}$ that form an antichain $\left \{x,y,z\right \} \subseteq \mathbf{M}_{n,n}$ with $x \in S$ and $y,z \in T$.
\end{definition}

Observe that by commutativity of lattice operations, $(x,y,z)$ is an $ST$-breaking antichain if and only if $(x,z,y)$ is.

\subsubsection{Theoretical Results 1: Generalization of 5 Types.}
We are now ready to prove that the constructions 2--5 in Remark \ref{pair-types} generalize to all $\mathbf{M}_{n,n}$. Note the central role that the MCCD and $ST$-breaking antichains play in the proofs and keep in mind our implicit assumption that $n \geq 3$. To avoid tedious repetition, we begin with the following remark that establishes some common ground for the subsequent propositions.

\begin{remark} (Common arguments for Proposition \ref{s-link}-\ref{s-level-minus-link})
The following apply in Propositions \ref{s-link}-\ref{s-level-minus-link} and their proofs:
\begin{enumerate}
\item $L_1, L_2$ denote the two levels of $\mathbf{M}_{n,n}$.
\item The $ST$-distributivity of the pair considered follows directly from the MCCD via a case by case analysis of the ways of choosing $s$, $t_1$, and $t_2$ from it. The verification for this is purely computational and is thus omitted. 
\item As a result of the previous item, the proof reduces to showing the maximality of the $ST$-pair in question.
\end{enumerate}
\end{remark}

\begin{proposition} (S-link). \label{s-link}
Let $l,x \in L_1$ where $l$ is the link. If $S = \left \{l\right \}$ and $T = L_2 \cup \left \{x\right \}$, then $\left(S,T\right)$ is a maximal $ST$-pair of $\mathbf{M}_{n,n}$.
\begin{proof}
Note that the only elements left in $\mathbf{M}_{n,n}^{*}$ are those in $L_1 \setminus \left \{l,x\right \}$. Let $y \in L_1 \setminus \left \{l,x\right \}$. If $y \in S$, then $\left ( y,x,z\right )$ is an $ST$-breaking antichain for any $z \in L_2$ except the link. Similarly, if $y \in T$, then $\left ( l,x,y\right )$ is an $ST$-breaking antichain. Hence, neither $S$ nor $T$ can be expanded.
\qed
\end{proof}
\end{proposition}

\begin{proposition} (S-2-links). \label{s-2-links}
Let $l_1$ be $L_1$'s link and $l_2$ be $L_2$'s link. Suppose $x \in L_1 \setminus \left \{l_1\right \}$ and $y \in L_2 \setminus \left \{l_2\right \}$. If $S = \left \{l_1, l_2\right \}$ and $T = \left \{x,y\right \}$, then $\left(S,T\right)$ is a maximal $ST$-pair of $\mathbf{M}_{n,n}$.
\begin{proof}
Note that the only elements left in $\mathbf{M}_{n,n}^{*}$ are those in 
\begin{equation}
\left(L_1 \setminus \left \{l_1,x\right \}\right) \cup \left(L_2 \setminus \left \{l_2,y\right \}\right). 
\end{equation}
Let $z$ be an element in this union. If $z \in S$, then $\left ( z,x,y\right )$ is an $ST$-breaking antichain. Now, suppose $z \in T$. Then either $\left ( l_1,x,z\right )$ or $\left ( l_2,y,z\right )$ is an $ST$-breaking antichain depending on whether $z \in L_1$ or $z \in L_2$.
\qed
\end{proof}
\end{proposition}

\begin{proposition} (S-level). \label{s-level}
Let $l,x \in L_1$ where $l$ is the link. If $S = L_2$ and $T = \left \{l,x\right \}$, then $\left(S,T\right)$ is a maximal $ST$-pair of $\mathbf{M}_{n,n}$.
\begin{proof}
Note that the only elements left in $\mathbf{M}_{n,n}^{*}$ are those in $L_1 \setminus \left \{l,x\right \}$. Suppose $y \in L_1 \setminus \left \{l,x\right \}$. If $y \in S$, then $\left ( y,x,l\right )$ is an $ST$-breaking antichain. If $y \in T$, then $\left ( z,x,y\right )$ is an $ST$-breaking antichain for all $z \in S = L_2$ except the link. 
\qed
\end{proof}
\end{proposition}

Observe that the only difference between $ST$-pairs of the type S-level-minus-link and those of type S-level is that pairs of type S-level-minus-link take the link element in the $S$ of a type S-level and add it to its $T$.

\begin{proposition} (S-level-minus-link). \label{s-level-minus-link}
Let $l_1, l_2$ be the respective links of $L_1, L_2$. Suppose $x \in L_1 \setminus \left \{l_1\right \}$. If $S = L_2 \setminus \left \{l_2\right \}$ and $T =  \left \{l_1, l_2, x\right \}$, then $\left(S,T\right)$ is a maximal $ST$-pair of $\mathbf{M}_{n,n}$. Dually, we can also suppose that $S = L_1 \setminus \left \{l_1\right \}$ and $T =  \left \{l_1, l_2, y\right \}$ for some $y \in L_2 \setminus \left \{l_2\right \}$ and the result will also hold. We omit it from the proof because it uses the same argument.
\begin{proof}
The same argument used in Proposition \ref{s-level} can be applied with $l = l_1$, adjusting only for the fact that $S$ no longer contains the link of $L_2$.
\qed
\end{proof}
\end{proposition}

\subsubsection{Theoretical Results 2: No More Pairs.}

Having established that all 5 types in Remark \ref{pair-types} generalize, we now show that we have identified all ways of constructing a maximal $ST$-pair of $\mathbf{M}_{n,n}$. The intuition behind this is that although the number of elements of $\mathbf{M}_{n,n}$ increases as $n$ does, many of the new triples are $ST$-breaking antichains since only the links, which are fixed at 2 for all $n$, may be comparable with the new elements. We begin with the following lemma, which greatly reduces the subsets of $\mathbf{M}_{n,n}^{*}$ to consider when choosing $T$. Its proof uses $\mathbb{N}_n$ to denote the set of the first $n$ natural numbers, $\{1, 2, ..., n\}$.

\begin{lemma} \itparen{Restriction on $T$} \label{t_restrict} 
If $S, T \subseteq \mathbf{M}_{n,n}^{*}$ are disjoint subsets such that $\mathbf{M}_{n,n}$ is $ST$-distributive, then $T$ cannot contain four elements, two from each level of $\mathbf{M}_{n,n}$.

\begin{proof}
We show this by contradiction. Suppose that $T \supseteq \{a_i, a_j, b_k, b_l \}$ for $i \neq j$ and $k \neq l$. If $a_m \in S$ for any $m \in \mathbb{N}_n \setminus \{i,j\}$, then $( a_m, a_i, a_j )$ is an $ST$-breaking antichain, contradicting the $ST$-distributivity of $\mathbf{M}_{n,n}$. A similar argument works for any $b_m \in S$ distinct from $b_k$ and $b_l$. This implies that $S = \emptyset$; contradiction. Therefore, $T \not\supseteq \{a_i, a_j, b_k, b_l \}$ for $i \neq j$ and $k \neq l$.
\end{proof}
\end{lemma}

Without further delay, we complete the characterization of maximal $ST$-pairs of $\mathbf{M}_{n,n}$ in Theorem \ref{conjecture_proof}. Table \ref{formal_pair_types} provides a formal description of the 5 types of maximal $ST$-pairs of $\mathbf{M}_{n,n}$ of Remark \ref{pair-types}, each with its possible variations. It will be a helpful reference when following the different cases of the proof of Theorem \ref{conjecture_proof}. 

\begin{table}
\begin{center}
\begin{tabular}{| c | c | c | c |}
\hline
Type 1 & \multicolumn{3}{c|}{T-chain} \\
\hline
\multirow{5}{*}{4 forms} & $S$ & $T$ & index restrictions  \\
\cline{2-4}
& $\{a_1, ..., a_{i-1}, a_{i+1}, ..., a_{n}, b_1, ... , b_n\}$ & $\{a_i\}$ & none \\
\cline{2-4}
& $\{a_1, ..., a_n, b_1, ..., b_{i-1}, b_{i+1}, ..., b_n\}$ & $\{b_i\}$ & none \\
\cline{2-4}
& $\{a_1, ..., a_{i-1}, a_{i+1}, ..., a_{n}, b_2, ... , b_n\}$ & $\{a_i, b_1\}$ & none \\
\cline{2-4}
& $\{a_1, ..., a_{n-1}, b_1, ..., b_{i-1}, b_{i+1}, ..., b_n\}$ & $\{a_n, b_i\}$ & none  \\
\hline
Type 2 & \multicolumn{3}{c|}{S-link} \\
\hline
\multirow{3}{*}{2 forms} & $S$ & $T$ & index restrictions \\
\cline{2-4}
& $\{b_1\}$ & $\{a_1, ..., a_n, b_i\}$ & $i \neq 1$ \\
\cline{2-4}
& $\{a_n\}$ & $\{a_i, b_1, ..., b_n\}$& $i \neq n$ \\
\hline
Type 3 & \multicolumn{3}{c|}{S-2-links} \\
\hline
\multirow{2}{*}{1 form} & $S$ & $T$ & index restrictions \\
\cline{2-4}
& $\{a_n,b_1\}$ & $\{a_i, b_j\}$ & $i \neq n, j \neq 1$ \\
\hline
Type 4 & \multicolumn{3}{c|}{S-level} \\
\hline
\multirow{3}{*}{2 forms} & $S$ & $T$ & index restrictions \\
\cline{2-4}
& $\{b_1, ..., b_n \}$ & $\{a_i, a_n \}$ &  $i \neq n$ \\
\cline{2-4}
& $\{a_1, ..., a_n \}$ & $\{b_1, b_i \}$ & $i \neq 1$ \\
\hline
Type 5 & \multicolumn{3}{c|}{S-level-minus-link} \\
\hline
\multirow{3}{*}{2 forms} & $S$ & $T$ & index restrictions\\
\cline{2-4}
& $\{b_2, ..., b_n \}$ & $\{a_i, a_n, b_1 \}$ & $i \neq n$ \\
\cline{2-4}
& $\{a_1, ..., a_{n-1} \}$ & $\{a_n, b_1, b_i \}$ & $i \neq 1$ \\
\hline
\end{tabular}
\caption{Formal descriptions of the 5 Types of Maximal $ST$-pairs of $\mathbf{M}_{n,n}$} \label{formal_pair_types}
\end{center}
\end{table}

\begin{theorem} \itparen{Characterization of maximal $ST$-pairs of $\mathbf{M}_{n,n}$} \label{conjecture_proof}
For all $n \in \mathbb{N}$ such that $n \geq 3$, the maximal $ST$-pairs of $\mathbf{M}_{n,n}$ all fall into one of the 5 types listed in Remark \ref{pair-types}  \itparen{or Table \ref{formal_pair_types}}. 

\begin{proof}
We show that if $(S,T)$ is a maximal $ST$-pair of $\mathbf{M}_{n,n}$, then $(S,T)$ is of one of Types 1-5. Suppose that $(S,T)$ is a a maximal $ST$-pair of $\mathbf{M}_{n,n}$. If $T$ is a chain, then $(S,T)$ is of Type 1 and we are done. For the remainder of the proof, we assume that $T$ is not a chain. By Lemma \ref{t_restrict}, $T$ cannot contain two elements of both levels of $\mathbf{M}_{n,n}$ at the same time. Hence, there are two cases:

\begin{enumerate}
\item $T$ is contained in one level: $T \subseteq \{a_1, ..., a_n \}$ or $T \subseteq \{b_1, ..., b_n \}$.
\item $T$ is contained in one level except for exactly one element: \\
$T \subseteq \{a_i, b_1, ..., b_n \}$ with $a_i \in T$ or $T \subseteq \{a_1, ..., a_n, b_i \}$ with $b_i \in T$.
\end{enumerate}

\textbf{\underline{Case 1}}: $T$ is contained in one level. Consider the case $T \subseteq \{a_1, ..., a_n \}$. The case $T \subseteq \{b_1, ..., b_n \}$ is done similarly. Since $T$ is not a chain, $|T| \geq 2$. That is $T$ contains $a_i$ and $a_j$ with $i \neq j$. Then no other $a_k$ can be in $S$ because $( a_k, a_i, a_j )$ would be an $ST$-breaking antichain. Thus, $S \subseteq \{b_1, ..., b_n\}$. This reduces the problem to determining how many $a_i$'s and $b_j$'s can be put in $T$ and $S$ respectively. \\

\indent To begin, we have $b_1 \in S$ by MCCD since $b_1 \geq a_i$ for all $i$. Now, observe that if $T$ contains two non-link and distinct $a_i$ and $a_j$, then $S = \{b_1\}$ because $( b_k, a_i, a_j )$ is an $ST$-breaking antichain for any $k \neq 1$.
This implies that $T = \{a_1, ..., a_n \}$ given the maximality of our pair and our assumption under Case 1. However, this pair is properly contained  in a pair of Type 2 contradicting its maximality. Hence, $T$ cannot have two distinct non-link elements under Case 1.\\

\indent We are left with the possibility that $|T| = 2$ with $a_n \in T$. In this case, we must  add all of the $b_i$'s to $S$ because $a_n \leq b_i$ for all $i$ \itparen{MCCD}. Thus, $S =\{b_1, ..., b_n\}$ and $T = \{a_i, a_n\}$ with $i \neq n \implies (S,T)$ is a maximal $ST$-pair of Type 4. \\

\indent \textbf{\underline{Case 2}}: $T$ is contained in one level except for exactly one element. There are two cases: \itparen{a} $|T| = 2$ and \itparen{b} $|T| > 2$. \\

\indent \emph{\underline{Sub-case 2a}}: If $|T| = 2$, then $T = \{a_i, b_j\}$. Since $T$ is not a chain, we must have $i \neq n$ and $j \neq 1$. Hence, $T$ is an antichain of two elements. This implies that for any non-link $a_k$ or $b_k$, $( a_k, a_i, b_j )$ and $( b_k, a_i, b_j )$ is an $ST$-breaking antichain. Thus, $S \subseteq \{a_n, b_1\}$. Given that $a_n \leq b_j$ and $b_1 \geq a_i$, the MCCD allows us to add both links to $S$. Thus we have $S = \{a_n, b_1\}$ and $T = \{a_i, b_j\}$ with $i \neq n$ and $j \neq 1 \implies (S,T)$ is a maximal $ST$-pair of Type 3. \\

\indent \emph{\underline{Sub-case 2b}}: Suppose $|T| > 2$. We work with the case that $T \subseteq \{a_i, b_1, ..., b_n \}$ with $a_i \in T$. The case $T \subseteq \{a_1, ..., a_n, b_i \}$ with $b_i \in T$ is done similarly. We have that T contains $\{a_i, b_j, b_k\}$ for $j \neq k$ because $|T| > 2$. Then $S$ cannot contain any $b_l$'s because $( b_l, b_j, b_k )$ is an $ST$-breaking antichain for any $b_l$ distinct from $b_j$ and $b_k$. Thus, 
\begin{equation} \label{s_in_a-level}
S \subseteq \{a_1, ..., a_{i-1}, a_{i+1}, ..., a_n \}.
\end{equation}
Again, we have two cases: \itparen{i} $a_i$ is a link $(i = n)$ and \itparen{ii} $a_i$ is a non-link $(i \neq n)$. \\

\indent \underline{Sub-sub-case 2b-i}: If $a_i$ is a link \itparen{$i = n$}, then $S \subseteq \{a_1, ..., a_{n-1} \}$ by \itparen{\ref{s_in_a-level}}. Note, however, that for any non-link $a_l$ \itparen{i.e. $l \neq n$}, $( a_l, b_j, b_k )$ is an $ST$-breaking antichain if both $b_j$ and $b_k$ are non-link elements. As a result, $T$ can have at most one non-link $b_j$ if we want $S$ to be non-empty. In addition,  $b_1$ can be added to $T$ because $b_1 \geq a_l$ \itparen{MCCD}. We can also add one non-link $b_k$ without affecting the content of $S$ since $a_n \leq b_k$ \itparen{MCCD}. Thus, $T = \{a_n, b_1, b_k \}$ with $k \neq 1$. Maximality then gives $S = \{a_1, ...,  a_{n-1}\}$. Therefore, $(S,T)$ is a maximal $ST$-pair of Type 5.\\

\indent \underline{Sub-sub-case 2b-ii}: If $a_i$ is a non-link \itparen{$i \neq n$}, then we know for sure that $T$ has two non-link elements because it will have at least one $b_j$ with $j \neq 1$ \itparen{recall $|T| > 2$} and $a_i$. Then for any non-link $a_l$ we attempt to place in $S$, $( a_l, a_i, b_j )$ is an $ST$-breaking antichain. Given \itparen{\ref{s_in_a-level}}, the only possible element of $S$ is $a_n$. Since $a_n \leq b_k$ for all $k$, MCCD implies that $S = \{a_n\}$. Maximality forces $T = \{a_i, b_1, ..., b_n \}$ with $i \neq n \implies (S,T)$ is a maximal $ST$-pair of Type 2. \\

\indent $\therefore (S,T)$ is of Type $m$ for some $m \in \{1,2,3,4,5\}$. 
\end{proof}

\end{theorem}

\subsubsection{Combinatorial Corollary.}

As a corollary of Theorem \ref{conjecture_proof}, we can count how many maximal $ST$-pairs $\mathbf{M}_{n,n}$ has for each $n$. Observe that the 5 types of pairs are all pairwise disjoint, e.g., it is not possible for a pair to be both a T-chain pair and an S-link pair. Therefore, we just have to count how many of each type there are using basic combinatorics and then add them up. We conclude this section with the proof of this count.

\begin{corollary} \itparen{Count of maximal $ST$-pairs of $\mathbf{M}_{n,n}$}
Let $n \geq 3$. Then $\mathbf{M}_{n,n}$ has $n^{2} + 8n - 6$ maximal $ST$-pairs which are divided among the 5 types as follows:
\begin{enumerate}
\item T-chain: $4n - 1$,
\item S-link: $2n - 2$,
\item S-2-links: $n^{2} - 2n + 1$,
\item S-level: $2n - 2$,
\item S-level-minus-link: $2n - 2$.
\end{enumerate}

\begin{proof}
We count how many pairs there are of each type and then sum them up. We point to Table  \ref{formal_pair_types} as a useful reference for following this proof. \\

\indent \underline{T-chain}: Since $S$ is completely determined by $T$ \itparen{$S = \mathbf{M}_{n,n}^{*} \setminus T$}, then there is exactly one T-chain pair for each chain of $\mathbf{M}_{n,n}^{*}$. There are $2n$ size 1 chains, $2n-1$ size 2 chains, and no larger chains. The number of size 1 chains is the number of elements of $\mathbf{M}_{n,n}^{*}$: $2n$. As for the size 2 chains, note that they are all of the following two forms: $\{a_n, b_i\}$ and $\{a_i,b_1\}$. Each of these forms has $n$ pairs, giving $2n$. However, we must subtract one because $\{a_n,b_1\}$ is counted twice. Adding these up gives $4n-1$ pairs of type T-chain. \\

\indent \underline{S-link}: Note that $S$ consists of only $a_n$ or $b_1$. Given $S$, $T$ automatically has all of the elements of the level not containing $S$'s link, but may differ by one non-link element of the level containing $S$'s link. Hence,
\begin{align}
\#pairs &= (\#\text{link choices})(\#\text{choices of non-link in $T$}) \\
&= 2 \times (n-1) \\
&= 2n - 2.
\end{align} 

\indent \underline{S-2-links}: Observe that $S$ is fixed as $\{a_n,b_1\}$. Thus, the number of pairs of this type are determined by the number of ways of choosing $T$. Since $T$ consists of one non-link $a_i$ and one non-link $b_j$, we have
\begin{align}
\#pairs &= (\#\text{non-link choices of $a_i$})(\#\text{non-link choices of $b_j$}) \\
&= (n-1) \times (n-1) \\
&= n^{2} - 2n + 1.
\end{align} 

\indent \underline{S-level}: Since $S$ is one of the two levels of $\mathbf{M}_{n,n}$, it has two options. Given $S$, $T$ has two elements from the other level. One of these must be the link while the other may be any of the non-links. As a result, 
\begin{align}
\#pairs &= (\#\text{levels})(\#\text{choices for $T$'s non-link}) \\
&= 2 \times (n-1) \\
&= 2n - 2.
\end{align} 

\indent \underline{S-level-minus-link}: Observe that each $ST$-pair of the type S-level-minus-link is obtained from a type $S$-level pair by taking the link element in the $S$ of a type S-level and adding it to its $T$. Therefore, there is a 1-1 correspondence between these pairs and hence $2n - 2$ pairs of type S-level-minus-link. \\

\indent \underline{Total}: It remains to add up the totals of pairs of each type. It can be verified that 
\begin{equation}
(4n-1) + (2n-2) + (n^{2}-2n+1) + (2n-2) + (2n-2)= n^{2} + 8n - 6,
\end{equation}
which is the total amount of pairs initially claimed.
\end{proof}
\end{corollary}

\subsection{Connection to Distributive Elements}  \label{distr-elem-connect}

Now that we are done with our main results regarding $ST$-distributive lattices, we relate said lattices with distributive elements of lattices. These have a history of changing in meaning. In Garret Birkhoff's {\it Lattice Theory} \cite{Bi}, they are defined as being synonymous to neutral elements (Chapter III, Section 9). However, distributive and neutral elements are more recently defined as two separate concepts, as can be seen in George Gr\"{a}tzer's {\it Lattice Theory: Foundation}  \cite{G2011} (Chapter III, Section 2). Nevertheless, we will define a distributive element as in \cite{Bi} for reasons that will be clarified shortly. Thus, a distributive element is an element of a lattice with the property that any three-element subset of the lattice containing it generates a distributive sublattice.

\begin{definition} \itparen{distributive element} \label{def_distr_elem}
A distributive element of a lattice $L$ is an element $a \in L$ such that for all $x,y \in L$, $[\{a,x,y\}]$ is distributive where $[\{a,x,y\}]$ is the sublattice of $L$ generated by $\{a,x,y\}$. The set of all distributive elements of $L$ is denoted $\Distr{L}$.
\end{definition}

\indent We now consider two weaker variations of distributive elements that we will call meet distributive elements and join distributive elements. These really are, respectively, dually distributive elements and distributive elements as defined in \cite{G2011} in disguise. We renamed them to make their connection to $ST$-distributivity more apparent. This is also why we define distributive elements as in \cite{Bi}.

\begin{definition} \itparen{meet distributive element} \label{def_mdistr_elem}
A meet distributive element of a lattice $L$ is an element $a \in L$ such that 
\begin{equation} \label{def_mdistr_elem_eq}
a \land (x \lor y) = (a \land x) \lor (a \land y)
\end{equation}
for all $x,y \in L$. The set of all meet distributive elements of $L$ is denoted $\mDistr{L}$.
\end{definition}

\begin{definition} \itparen{join distributive element} \label{def_jdistr_elem}
A join distributive element of a lattice $L$ is an element $a \in L$ such that 
\begin{equation} \label{def_jdistr_elem_eq}
a \lor (x \land y) = (a \lor x) \land (a \lor y)
\end{equation}
for all $x,y \in L$. The set of all join distributive elements of $L$ is denoted $\jDistr{L}$.
\end{definition}

\indent Note that all elements of a distributive lattice are trivially distributive, meet distributive, and join distributive. In addition, it is immediate from Definitions \ref{def_distr_elem}-\ref{def_jdistr_elem} that all distributive elements are meet distributive and join distributive. We use this fact to  give the first non-trivial example of these three types of elements.

\begin{example} \itparen{lattice identities}
The identity elements $0$ and $1$ of any lattice $L$ are distributive. This is because $[\{0,x,y\}]$ is either a chain, $\mathbf{M}_2$, or $\mathbf{1} \oplus \mathbf{M}_2$ with dual reasoning for $[\{1,x,y\}]$. This results in $0$ and $1$ both being meet and join distributive as well.
\end{example}

\indent Next, we make a few additional comments regarding the relationships between these three types of lattice elements. It follows naturally from their definitions that meet distributive elements and join distributive elements are dual concepts. However, they are not equivalent. Furthermore, although all distributive elements are both meet distributive and join distributive, it is possible for an element to be both meet distributive and join distributive but not distributive. \\

\indent  We now give three more examples, two of which will illustrate the relations just discussed.
They involve the lattices $\mathbf{M}_3$ and $\mathbf{N}_5$, which we repeat in Figures \ref{m3_distr_elem} and \ref{n5_distr_elem} for ease of reference. We omit the computations involved because they are elementary.

\begin{example} \itparen{meet distributive $\centernot\iff$ join distributive}
An element that is meet distributive but not join distributive is $v \in \mathbf{N}_5$. By  Property \ref{eff_crit} and idempotency of lattice operations, proving that $v$ is meet distributive reduces to checking Equation \itparen{\ref{def_mdistr_elem_eq}} for the case $a = v$, $x = u$, and $y = w$. On the other hand, to see that it is not join distributive, apply the same variable assignments in Equation \itparen{\ref{def_jdistr_elem_eq}}. It can be similarly shown that $u \in \mathbf{N}_5$ is join distributive but not meet distributive.
\end{example}

\begin{example} \itparen{meet and join distributive but not distributive}
The element $w$ of $\mathbf{N}_5$ is both meet distributive and join distributive but not distributive. To see that it is both meet and join distributive it suffices to consider only the cases where $x$ and $y$ are both distinct from each other and from $w$ in Equations \itparen{\ref{def_mdistr_elem_eq}} and \itparen{\ref{def_jdistr_elem_eq}}. To see that it is not distributive note that $[\{w,u,v\}] = \mathbf{N}_5$. 
\end{example}

\begin{example} \itparen{complete non-example}
The element $a \in \mathbf{M}_3$ is neither meet distributive nor join distributive. To see this, it is enough to substitute $a = a$, $x = b$, and $y = c$ in Equations \itparen{\ref{def_mdistr_elem_eq}} and \itparen{\ref{def_jdistr_elem_eq}}. This also implies that $a$ is not distributive, which can also be verified  by observing that $[\{a,b,c\}] = \mathbf{M}_3$.
\end{example}

\begin{figure}
\begin{minipage}{0.45\textwidth}
\centering
\begin{tikzpicture}
[roundnode/.style={circle, draw=black, fill=white, very thin, text width =4.00mm, inner sep = 0pt, text centered}]
 \node [roundnode] (1) at (0,2) {$1$};
  \node [roundnode] (a) at (-1,1) {$a$};
  \node [roundnode] (b) at (0,1) {$b$};
  \node [roundnode] (c) at (1,1) {$c$};
  \node [roundnode] (0) at (0,0) {$0$};
 \draw (1) -- (a) -- (0) -- (b) -- (1) -- (c) -- (0);
\end{tikzpicture}
\caption{The diamond $\mathbf{M}_3$} \label{m3_distr_elem}
\end{minipage} \hfill
\begin{minipage}{0.45\textwidth}
\centering
\begin{tikzpicture}
[roundnode/.style={circle, draw=black, fill=white, very thin, text width =4.00mm, inner sep = 0pt, text centered}]
 \node [roundnode] (1) at (0,2) {$1$};
  \node [roundnode] (u) at (-1,1.4) {$u$};
  \node [roundnode] (v) at (-1,.6) {$v$};
  \node [roundnode] (w) at (1,1) {$w$};
  \node [roundnode] (0) at (0,0) {$0$};
 \draw (1) -- (u) -- (v) -- (0) -- (w) -- (1);
\end{tikzpicture}
\caption{The pentagon $\mathbf{N}_5$} \label{n5_distr_elem}
\end{minipage}
\end{figure}

\indent Finally, we mention the connection between distributive elements and $ST$-distributive lattices. This is that Definitions \ref{def_distr_elem}-\ref{def_jdistr_elem} provide choices of $S \subseteq L$ that guarantee $ST$-distributivity for $T = L$. In particular, we have that any lattice $L$ is 
\begin{enumerate}
\item $(\mDistr{L})L$-meet distributive,
\item $(\jDistr{L})L$-join distributive,
\item $(\Distr{L})L$-distributive.
\end{enumerate}
%

\subsection{Connection to SS-lattices} \label{ss-lat-connect}

Before moving on to $ST$-modular lattices, we place our $ST$-distributive lattices in the context of Stanley's SS-lattices (see \cite{S}). Our main point is that SS-lattices provide examples of $ST$-distributive lattices but $ST$-distributive lattices are more general. First, we recall the definition of an SS-lattice.

\begin{definition} (SS-lattice \cite{S}).
Consider a lattice $L$ with a maximal chain $\Delta$. $\left(L,\Delta\right)$ is called a supersolvable lattice (SS-lattice) if for all chains $K$ of $L$, $\left[\Delta,K\right]$ is distributive. Here, $\left[\Delta,K\right]$ denotes the sublattice of $L$ generated by $\Delta \cup K$.
\end{definition}

From this definition, we can observe how we can obtain examples of $ST$-distributive lattices from it: if $\left(L,\Delta\right)$ is an SS-lattice, then $L$ is $\Delta K$-distributive and $K \Delta$-distributive for all its chains $K$. On the other hand, $ST$-distributive lattices are a bigger class because the conditions specified on Definition \ref{st-distr} are weaker than those required for an SS-lattice. First, we let $S$ and $T$ be any subsets of the initial lattice, not specifically chains. In addition, we only ask that $S$ distributes into $T$ rather than demanding that the generated sublattice $\left[S,T\right]$ be distributive. This concludes our treatment of $ST$-distributive lattices.

\section{$ST$-modular Lattices} \label{sec_st-mod}

This section explores $ST$-modular lattices. As we shall see, these expand the class of $ST$-distributive lattices. There are two parts to this section. First, Subsect. \ref{basic-def-st-mod} gives the basic definitions and some examples. Afterwards, Subsect. \ref{conv_app} suggests a potential application to convex sets. In particular, we will use $ST$-modularity to reprove some classical results on convex polytopes. We remark that this connection to convexity is the original motivation for defining $ST$-modular lattices, from which $ST$-distributive lattices follow by analogy.

\subsection{Basic Definitions} \label{basic-def-st-mod}

We now propose a generalization of modular lattices that is analogous to $ST$-distributive lattices: $ST$-modular lattices. We first define them, then explain how they relate to $ST$-distributive lattices, and conclude with some illustrative examples. The reader should observe the parallels between this subsection and Subsect. \ref{basic-def-st-distr}: the first being the similarity between Definition \ref{st-distr} and the following definition:


\begin{definition} ($ST$-modular Lattice).
Given a lattice $L$ with $S, T \subseteq L$, we define:

\begin{itemize}
\item \underline{$ST$-meet Modular Lattice}: $L$ is said to be $ST$-meet modular if for all  $s \in S$ and $t_1, t_2 \in T$,
\begin{equation}
s \geq t_2 \implies s \land \left(t_1 \lor t_2\right) = \left(s \land t_1\right) \lor t_2 \; .
\end{equation}

\item \underline{$ST$-join Modular Lattice}: $L$ is said to be $ST$-join modular if for all  $s \in S$ and $t_1, t_2 \in T$
\begin{equation}
s \leq t_2 \implies s \lor \left(t_1 \land t_2\right) = \left(s \lor t_1\right) \land t_2 \; .
\end{equation}

\item \underline{$ST$-modular Lattice}: $L$ is said to be $ST$-modular if it is both $ST$-meet modular and $ST$-join modular.
\end{itemize}
\end{definition}

\begin{figure}
\centering
\begin{tikzpicture}
[roundnode/.style={circle, draw=black, fill=white, very thin, text width =4.00mm, inner sep = 0pt, text centered}]
 \node [roundnode] (1) at (0,3.5) {$1$};
  \node [roundnode] (a) at (0,2) {$a$};
  \node [roundnode] (b) at (-2,2) {$b$};
  \node [roundnode] (c) at (2,2) {$c$};
  \node [roundnode] (d) at (-1,1) {$d$};
  \node [roundnode] (e) at (1,1) {$e$};
  \node [roundnode] (0) at (0,0) {$0$};
 \draw (1) -- (b) -- (d) -- (0) -- (e) -- (c) -- (1) -- (a) -- (d);
 \draw (a) -- (e);
\end{tikzpicture}
\caption{Lattice $L$ of Examples \ref{st-mod-not-st-distr}, \ref{mod_st-meet-neq-st-join}, and \ref{mod_st-neq-ts}} \label{counter-ex-lattice2}
\end{figure}
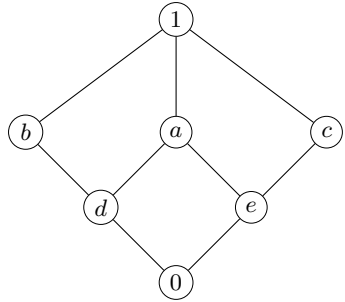


We now briefly discuss $ST$-modular lattices. For starters, note how these generalize modular lattices: $ST$-modular lattices for $S = T = L$. Next, before giving some concrete examples, we examine the relationship between $ST$-modular and $ST$-distributive lattices. We establish that $ST$-distributive lattices are $ST$-modular. This follows directly from the fact that distributive lattices are modular \cite{Bi,DP,G2011}. 

\begin{proposition} ($ST$-distributive $\implies$ $ST$-modular). \label{st-distr-implies-mod}
Let $L$ be a lattice with $S,T \subseteq L$. If $L$ is $ST$-distributive, then $L$ is $ST$-modular

\begin{proof}
Suppose $L$ is $ST$-distributive for $S,T \subseteq L$ with $s \in S$ and $t_1, t_2 \in T$. If $s \geq t_2$, then the $ST$-distributivity of $L$ and the Connecting Lemma \cite{Bi,DP} imply

\begin{equation}
s \land \left(t_1 \lor t_2\right) = \left(s \land t_1\right) \lor \left(s \land t_2\right) = \left(s \land t_1\right) \lor t_2 \; .
\end{equation}

Similarly, if $s \leq t_2$,

\begin{equation}
s \lor \left(t_1 \land t_2\right) = \left(s \lor t_1\right) \land \left(s \lor t_2\right) = \left(s \lor t_1\right) \land t_2 \; .
\end{equation} 
\qed
\end{proof}
\end{proposition}

Thus, we have that the $ST$-distributive lattice from Example \ref{st-distr-ex} is also an $ST$-modular lattice. Another $ST$-modular lattice is given in Example \ref{st-mod-not-st-distr}. It is also a counter-example to the converse of Proposition \ref{st-distr-implies-mod}. Hence, the class of $ST$-modular lattices is greater than that of $ST$-distributive ones. This is consistent with prior knowledge that not all modular lattices are distributive  \cite{Bi,DP,G2011}.

\begin{example} ($ST$-modular $\centernot\implies ST$-distributive). \label{st-mod-not-st-distr}
Recall the lattice $L$ from Fig. \ref{counter-ex-lattice}, repeated in Fig. \ref{counter-ex-lattice2} for convenience. Consider the subsets $S = \left \{b,d\right \}$ and $T =\left \{a,c\right \}$. We claim the $L$ is $ST$-modular but not $ST$-distributive.
\begin{proof}
\noindent \underline{$ST$-modular}: Note that there is no $s \in S$ and $t \in T$ such that $s \geq t$. Hence, $L$ is vacuously $ST$-meet modular. To show that it is $ST$-join modular, we only need to consider when $d \leq a$. This leads to 2 cases where the necessary equalities hold. Note the use of absorption in the rightmost equality of the first case.
\begin{eqnarray}
d \lor \left(a \land a\right) = d \lor a = a  = \left(d \lor a\right) \land a \qquad \: \\
d \lor \left(c \land a\right) = d \lor e = a =  1 \land a =  \left(d \lor c\right) \land a \; .
\end{eqnarray}

\underline{Not $ST$-distributive}: Note that $b \in S$ and $a,c \in T$ yet,
\begin{equation}
b \land \left(a \lor c\right) = b \land 1 = b \neq  d = d \lor 0 = \left(b \land a\right) \lor \left(b \land c\right)\; .
\end{equation}
\qed
\end{proof}
\end{example}

We finish this subsection with two more examples. These illustrate, respectively, that $ST$-modular lattices behave like $ST$-distributive lattices when it comes to (a) the distinction between their $ST$-meet and $ST$-join sub-definitions and (b) the ordered nature of the pair of subsets $\left(S,T\right)$.


\begin{example} ($ST$-meet Modular $\centernot\iff$ $ST$-join Modular). \label{mod_st-meet-neq-st-join}
Examine again the lattice $L$ from Fig. \ref{counter-ex-lattice2}. Consider the subsets $S = \left \{0,b\right \}$ and $T =\left \{c,d\right \}$. We claim the $L$ is $ST$-join modular but not $ST$-meet modular.
\begin{proof}
\noindent \underline{$ST$-join Modular}: Note that $0 \leq c$ and $0 \leq d$. This leads to 4 cases that must be verified. It can be seen from Table \ref{non_equiv_st-join_mod} that all of them the satisfy the modular law. \\

\begin{table}
\begin{center}
\caption{Computations to show the $L$ is $ST$-join modular in Example \ref{mod_st-meet-neq-st-join}} \label{non_equiv_st-join_mod}

\begin{tabular}{ c m{7mm} c }
\hline 
$0 \lor \left(c \land c\right) = c$ & & $0 \lor \left(c \land d\right) = 0$ \\
$\left(0 \lor c\right) \land c = c$ & & $\left(0 \lor c\right) \land d = 0$ \\
\hline
$0 \lor \left(d \land c\right) = 0$ & & $0 \lor \left(d \land d\right) = d$ \\
$\left(0 \lor d\right) \land c= 0$ & & $\left(0 \lor d\right) \land d = d$\\
\hline
\end{tabular} 
\end{center}
\end{table}

\underline{Not $ST$-meet Modular}: Observe that although $b \in S$, $c, d \in T$ and $b \geq d$, we have:
\begin{equation}
b \land \left(c \lor d\right) = b \land 1 = b \neq d = 0 \lor d = \left(b \land c\right) \lor d \; .
\end{equation}
\qed
\end{proof}
\end{example}


\begin{example} ($ST$-modular $\centernot\iff$ $TS$-modular). \label{mod_st-neq-ts}
We use the lattice $L$ from Fig. \ref{counter-ex-lattice2} one more time. Consider the subsets $S = \left \{b,c\right \}$ and $T =\left \{a,d\right \}$. We claim the $L$ is $ST$-modular but not $TS$-modular.
\begin{proof}
\noindent \underline{$ST$-modular}: Note that there is no $s \in S$ and $t \in T$ such that $s \leq t$. Hence, $L$ is vacuously $ST$-join modular. To show that it is $ST$-meet modular, we only need to consider when $b \geq d$. This leads to 2 cases where the necessary equalities hold. Note the use of absorption in the rightmost equality of the first case.
\begin{eqnarray}
b \land \left(d \lor d\right) = b \land d = d = \left(b \land d\right) \lor d \qquad \:\\
b \land \left(a \lor d\right) = b \land a = d = d \lor d = \left(b \land a\right) \lor d \; .
\end{eqnarray}

\underline{Not $TS$-modular}: To see this, note that although $d \in T$, $b,c \in S$, and $d \leq b$, we have 
\begin{equation}
d \lor \left(c \land b\right) = d \lor 0 = d \neq b = 1 \land b = \left(d \lor c\right) \land b \; .
\end{equation} \qed
\end{proof}
\end{example}

\subsection{Application to Convex Sets} \label{conv_app}

Having defined what is an $ST$-modular lattice, we will take a different route from that taken with $ST$-distributive lattices. Rather than searching for $ST$-modular lattices, we will propose an application of this relative modular property to convex sets. This application is the original motivation for this concept. For it we will need two ingredients: the  lattice of convex sets and a property of convex hulls in Proposition \ref{conv_hull_ident}. In the following, $\mathrm{conv}\left(A,B\right)$ denotes the convex hull of $A \cup B$ for two sets $A$ and $B$.


\begin{definition} (Lattice of Convex Sets \cite{B}).
The lattice of convex sets of $\mathbb{R}^{d}$ is the lattice $\langle 
L, \lor, \land \rangle$ where $L$ is the set of all convex sets of $\mathbb{R}^{d}$ with lattice operations
\begin{equation}
A \land B = A \cap B \qquad  A \lor B = \mathrm{conv}\left(A,B\right) \; .
\end{equation}
\end{definition}

We mention two facts of this lattice: that it is ordered by inclusion ($\subseteq$) and that it can be shown to be non-modular. The latter will make our proposed relative modularity non-trivial. Now, we present our second ingredient: a statement regarding the intersection of an affine set with the convex hull of two convex sets.


\begin{proposition} (Interesection of Affine Set and Convex Hull \cite{E}). \label{conv_hull_ident}
Suppose we have $A, B, C \subseteq \mathbb{R}^{d}$ with $A$ affine, $B$ and $C$ convex, and $A \supseteq C$. Then:
\begin{equation}
A \cap \mathrm{conv}\left(B,C\right) = \mathrm{conv}\left(A \cap B, C\right) \; .
\end{equation}
\begin{proof}
The inclusion $\supseteq$ is trivial because $A \cap \mathrm{conv}\left(B,C\right)$ is convex and contains both $A \cap B$ and $C = A \cap C$. For the other inclusion ($\subseteq$), let $x \in A \cap \mathrm{conv}\left(B,C\right)$. Then $x \in \mathrm{conv}\left(B,C\right)$ and since $\mathrm{conv}\left(B,C\right)$ is the union of all line segments connecting a point in $B$ with a point in $C$, there exist $b \in B$, $c \in C$ and $0 \leq \alpha \leq 1$ such that
\begin{equation} \label{line_expr}
x = \left(1 - \alpha\right) b + \alpha c \; .
\end{equation}

If $\alpha = 0$, then $x = b \in B$ and given that $x \in A$, $x \in A \cap B$ and we are done. Similarly, if $\alpha = 1$, then $x = c \in C$ and we are done. Otherwise, solve Eq. (\ref{line_expr}) for $b$ to get
\begin{equation}
b = \frac{1}{1-\alpha}x - \frac{\alpha}{1-\alpha}c \; .
\end{equation}

It follows that $b \in A$. This due to the fact that $x,c \in A$ (recall $C \subseteq A$), the sum of their coefficients is 1, and $A$ is affine. This means that the convex combination presentation in Eq. (\ref{line_expr}) is a convex combination of a point in $A \cap B$ and a point in $C$. Therefore, $x \in \mathrm{conv}\left(A \cap B, C\right)$.\qed
\end{proof}
\end{proposition}


We are ready to present our application of $ST$-modularity to convex sets. The main idea is that we can use our relative modularity to describe certain behaviors of convex sets. For instance, we can translate Proposition \ref{conv_hull_ident} into the language of the lattice of convex sets to get $ST$-meet modularity in said lattice.

\begin{corollary} ($ST$-meet Modularity in the Lattice of Convex Sets).
Let $L$ be the lattice of convex sets of $\mathbb{R}^{d}$. If $S $ is the set of all affine sets of $\mathbb{R}^{d}$ and $T = L$, then Proposition \ref{conv_hull_ident} simply says that $L$ is $ST$-meet modular. This is because for $A \in S$, $B, C \in T$ with $A \geq C$, it implies
\begin{equation}
A \land \left(B \lor C\right) = \left(A \land B\right) \lor C \; .
\end{equation}
\end{corollary}


The relevance of this connection is that this $ST$-meet modularity can be used to prove classical theorems on convex polytopes. We will provide two examples. Before doing so, we clarify some definitions and notations that will be used in the proofs. We define a \emph{convex polytope} as the convex hull of a finite set of points in $\mathbb{R}^{d}$ and a \emph{face} of a convex polytope as the intersection of the polytope with a supporting hyperplane (following \cite{E}). A \emph{vertex} is a 0-dimensional face. Given a hyperplane $H$ in $\mathbb{R}^{d}$, $H^{+}$ denotes one of the closed halfspaces of $\mathbb{R}^{d}$ that it defines and $int\left(H^{+}\right)$ its corresponding open halfspace. Finally, for a finite set of points $\left \{x_1, ..., x_n\right \}$, we write $\mathrm{conv}\left(x_1, ..., x_n\right)$ as short hand for $\mathrm{conv}\left(\left \{x_1, ..., x_n\right \}\right)$.

\begin{theorem} (Faces of Convex Polytope \cite{E}).
\begin{enumerate}
\item A convex polytope has a finite number of faces.
\item Each face is a convex polytope. 
\end{enumerate}
\newpage
\begin{proof}
We prove Statement 2 first and then derive Statement 1 from it. \\

\underline{Statement 2}: Let $P = \mathrm{conv}\left(x_1, ..., x_n\right)$ be a convex polytope and let $F$ be a face of $P$. Then there exists a supporting hyperplane $H$ of $P$ such that $F = P \cap H$ with $x_1, ..., x_j \in F$ and $x_{j+1}, ..., x_n \in int\left(H^{+}\right)$ for some $j \in \left \{1, ..., n \right \}$. Then
\begin{eqnarray}
&F &= H \cap P \\
& &= H \cap \mathrm{conv}\left(x_1, ..., x_n\right) \\
& &= H \cap \mathrm{conv}\left( \mathrm{conv}\left(x_{j+1}, ..., x_n\right), \mathrm{conv}\left(x_1, ..., x_j\right) \right) \; .
\end{eqnarray}
\indent Since $H$ is affine,  $\mathrm{conv}\left(x_1, ..., x_j\right)$ and $\mathrm{conv}\left(x_{j+1}, ... x_n\right)$ are both convex, and $H \supseteq \mathrm{conv}\left(x_1, ..., x_j\right)$, we can rewrite the right-side expression using Proposition \ref{conv_hull_ident}:
\begin{eqnarray}
&F &= \mathrm{conv}\left( H \cap \mathrm{conv}\left(x_{j+1}, ..., x_n\right), \mathrm{conv}\left(x_1, ..., x_j\right) \right) \\
& &= \mathrm{conv}\left( \emptyset, \mathrm{conv}\left(x_1, ..., x_j\right) \right) \\
& &= \mathrm{conv}\left(x_1, ..., x_j\right) \; .
\end{eqnarray} 
\indent Therefore, $F$ is a convex polytope. \\

\underline{Statement 1}: Now, observe that we just showed that each face of $P$ is the convex hull of a subset of $\left \{x_1, ..., x_n\right \}$, which is a finite set. Thus, the number of faces of $P$ must be finite.
\qed
\end{proof}
\end{theorem}

\begin{theorem} (Vertices of Convex Polytope \cite{E}).
A convex polytope is the convex hull of its vertices.
\begin{proof}
Let $P = \mathrm{conv}\left(x_1, ..., x_n\right)$ be a convex polytope. Suppose that this convex hull presentation of $P$ is minimal in the sense that 
\begin{equation}
P \neq \mathrm{conv}\left(x_1, ..., x_{i-1}, x_{i+1}, ..., x_n\right) = S_i
\end{equation}
for all $i$ in $1,..., n$. We show that this implies that each $x_i$ is a vertex and hence, $P$ is the convex hull of its vertices. \\

\indent First, note that by supposition, $x_i \not\in S_i$. Second, we have that $S_i$ is compact because it is the convex hull of a compact set: $\left \{x_1 ,..., x_{i-1}, x_{i+1}, ..., x_n\right \}$ Putting these together we get that $S_i$ is a convex and compact set not containing $x_i$. This implies the existence of a hyperplane $H$ passing through $x_i$ such that $S_i \subseteq int \left(H^{+}\right)$. We claim that $H$ is an appropriate hyperplane to establish that $x_i$ is a vertex.\\

\indent We must show that (a)  $\left \{x_i\right \} = H \cap P$ and (b) $P \subseteq H^{+}$. For the former, observe that
\begin{equation}
H \cap P = H \cap \mathrm{conv}\left(S_i,\left \{x_i\right \}\right)  \; .
\end{equation}
\indent Since $H$ is affine, $\left \{x_i\right \}$ and $S_i$ are both convex, and $H \supseteq \left \{x_i\right \}$, we can rewrite the right-hand side with Proposition \ref{conv_hull_ident} to obtain:
\begin{eqnarray}
&H \cap P &= \mathrm{conv}\left(H \cap S_i, \left \{x_i\right \}\right) \\
& &= \mathrm{conv}\left(\emptyset,\left \{x_i\right \} \right) \\
& &= \left \{x_i\right \} \; .
\end{eqnarray}

\indent It remains to show that $P$ is contained in $H^{+}$. Note that $\left \{x_i\right \}$ and $S_i$ are both subsets of $H^{+}$ and that $P = \mathrm{conv}\left(\left \{x_i\right \},S_i\right)$. Then $P \subseteq H^{+}$ because halfspaces are convex. Therefore, $x_i$ is a vertex of $P$.
\qed
\end{proof}
\end{theorem}

The main takeaway from all this is that $ST$-modularity may have applications to the study of convex sets and polytopes.

\section{Conclusion} \label{sec_concl}

This paper has introduced two new classes of lattices: $ST$-distributive lattices and $ST$-modular lattices. They extend distributive and modular lattices with corresponding relative properties satisfied by certain pairs of subsets of a lattice. \\

\indent Our main result is the full characterization of the maximal $ST$-pairs (of disjoint and non-empty subsets that exclude the identity elements) of the lattice family $\mathbf{M}_{n,n}$ for $n \geq 3$. To this end, we developed an algorithm that finds all such pairs and applied it $\mathbf{M}_{3,3}$ and $\mathbf{M}_{4,4}$. We found 27 and 42 maximal $ST$-pairs respectively. We classified them into 5 types of pairs based on how they are constructed. We then showed that these constructions generalize to all $\mathbf{M}_{n,n}$ for $n \geq 3$ (one of them even for all lattices) and that these are the only constructions possible in this family. Finally, we counted how many pairs of each type $\mathbf{M}_{n,n}$ has as a function of $n$. \\

\indent We also established some connection between our proposed classes of lattices and existing concepts in mathematics. We mentioned how to get $ST$-distributivity from distributive elements (as defined in \cite{Bi}) and renamed two subclasses of these meet and join distributive elements, respectively, to fit our new definitions.  We also discussed how our $ST$-distributive lattices include and extend Stanley's $SS$-lattices. Furthermore, we illustrated how $ST$-modularity can be used to describe properties of convex sets which can then be used to prove results on convex polytopes. \\

\indent Needless to say, there is plenty to be done in this new and exciting area of research. Naturally, we can continue our search for $ST$-distributive lattices with other families of lattices and also search for $ST$-modular ones. Of particular interest are cubical lattices (face lattices of hypercubes) and lattices of cut-complexes of hypercubes. We can also consider restricting our search to pairs of sublattices. All this can eventually lead to the development of a general theory of $ST$-distributive and $ST$-modular lattices with potential connections to other areas of mathematics (e.g. convex sets). The possibilities are endless!

\subsubsection*{Acknowledgements.} We would like to thank the Math Department at Florida Atlantic University and the organizers of SEICCGTC 2021 for all their effort in making the first virtual SEICCGTC possible.


%
%
\bibliographystyle{spmpsci}
\bibliography{reference}

\end{document}